\documentclass[12pt]{amsart}

\usepackage{amssymb}

\usepackage{graphicx}

\makeatletter
\@namedef{subjclassname@2010}{%
  \textup{2010} Mathematics Subject Classification}
\makeatother

\usepackage[T1]{fontenc}

\newtheorem{theorem}{Theorem}[section]

\theoremstyle{definition}
\newtheorem{definition}[theorem]{Definition}

\numberwithin{equation}{section}

\begin{document}


\baselineskip=17pt


\title[Supercommutator algebras of right (Hom-)alternative superalgebras]{Supercommutator (Hom-)superalgebras of Right 
(Hom-)Alternative Superalgebras}


\author[A. Nourou Issa]{A. Nourou Issa} 
\address{D\'epartement de Math\'ematiques, \\ Universit\'e d'Abomey-Calavi, \\ 01 BP 4521 Cotonou 01, BENIN}
\email{woraniss@yahoo.fr}

\date{}


\begin{abstract}

It is shown that the supercommutator superalgebra of a right alternative superalgebra is a Bol superalgebra. Hom-Bol 
superalgebras are defined and it is shown that they are closed under even self-morphisms. Any Bol superalgebra along 
with any even self-morphism is twisted into a Hom-Bol superalgebra. The supercommutator superalgebra of a right 
Hom-alternative superalgebra has a natural Hom-Bol structure. In order to prove this last result, 
the Hom-Jordan-admissibility of right Hom-alternative superalgebras is investigated and next Hom-Jordan supertriple 
systems are defined and their connection with Hom-Jordan superalgebras and Hom-Lie supertriple systems is considered.\\

\end{abstract}




\maketitle

\section{Introduction}

A {\it right alternative algebra} is an algebra satisfying the {\it right alternative identity} $(xy)y = x(yy)$. These algebras
were first considered in \cite{Alb}. For further studies on right alternative algebras, one may refer to \cite{Kl, Sk, 
Thed}.
\par
It turns out that right alternative algebras have close relations with a type of binary-ternary algebras called {\it Bol algebras}
which were introduced in \cite{Mik1} (see also references therein). In fact, it is proved \cite{Mik2} that any right alternative
algebra has a natural Bol algebra structure.
\par
The general theory of superalgebras started with the introduction of $\mathbb{Z}_2$-graded Lie algebras (i.e. Lie superalgebras)
coming from physics (see \cite{Kac1, Ree, Sch} and references therein for basics on Lie superalgebras). The important
role played by superalgebras in physics is beyond doubt. The $\mathbb{Z}_2$-graded
generalization of algebras is first extended to Jordan algebras in \cite{Kac2}. Next, alternative superalgebras were introduced in
\cite{ZS} whereas Maltsev superalgebras were introduced in \cite{Sh1}. A $\mathbb{Z}_2$-graded generalization of Bol algebras is 
considered in \cite{Ruk}.
\par
With the introduction of Hom-Lie algebras (see \cite{HLS, LS1, LS2}) began studies of Hom-type generalizations of usual
algebras. Apart from Hom-Lie algebras, first Hom-type algebras were defined in \cite{MS} while Hom-alternative and Hom-Jordan
algebras were defined in \cite{Ma} (see also \cite{Yau4} where Hom-Maltsev algebras were defined) and the Hom-type generalization of right
alternative algebras was considered in \cite{Yau2}. Other aspects of the theory of Hom-algebras are investigated in 
\cite{BM, MCL, Sheng, YSL, ZL}. It should be observed that, in general, the twisting map in a Hom-algebra is neither
injective nor surjective and when the twisting map is the identity map, then one recovers the ordinary (untwisted) algebraic structure. 
So ordinary algebras are viewed as Hom-algebras with the identity map as twisting map. Moving further in the theory of Hom-algebras, 
the Yau twisting principle of algebras is extended to binary-ternary algebras in \cite{Iss1} and next Hom-Bol 
algebras were defined in \cite{AI1}. As in the case of right alternative algebras, it is shown in \cite{AI2} that a 
Hom-Bol algebra structure can be defined on any multiplicative right Hom-alternative algebra.
\par
In this paper we extend, but with a different approach, the result in \cite{Mik2} to the cases of right alternative superalgebras and 
right Hom-alternative superalgebras.
Some basics on superalgebras are reminded in section 2 and some preliminary results are proved. In section 3 it is 
proved that any right alternative superalgebra has a natural Bol structure. The Hom-Jordan-admissibility
of multiplicative right Hom-alternative superalgebras is investigated in section 4.
In section 5, in order to deal with the Hom-version of the results from section 3, we first define Hom-Bol
superalgebras  and mention some construction results. In section 6 we extend to the $\mathbb{Z}_2$-graded Hom-type
setting the usual construction of Jordan triple systems from Jordan algebras by defining Hom-Jordan supertriple 
systems (the ungraded Hom-version of this construction
is considered in \cite{Yau3}). An example of a Jordan supertriple system is given. In section 7, relying on results 
from sections 4, 5, and 6, we finally prove that right Hom-alternative superalgebras are in fact Hom-Bol superalgebras by
showing first that any multiplicative Hom-Jordan superalgebra has a Hom-Jordan supertriple system structure and next
that any multiplicative Hom-Jordan supertriple system is a Hom-Lie supertriple system.
\par
All vector spaces and algebras are considered over a fixed ground field of characteristic not $2$ or $3$.

\section{Preliminaries}

A {\it superspace} (or a $\mathbb{Z}_2$-{\it graded space}) $V$ is  a direct sum $V = V_0\oplus V_1$, where $V_i$ are vector spaces. An element 
$x \in V_i$ ($i \in \mathbb{Z}_2$) is said to be {\it homogeneous of degree} $i$ and the degree of $x$ will be denoted by $\bar{x}$. 
\begin{definition}
 (i) Let $f: A \rightarrow A'$ be a linear map, where $A = A_0 \oplus A_1$ and $A' = A'_0 \oplus A'_1$ are superspaces. The
map $f$ is said to be {\it even} (resp. {\it odd}) if $f(A_i) \subset A'_i$ (resp. $f (A_i ) \subset A'_{i+1})$ for $i = 0,1$. 
\par
(ii) A ({\it multiplicative}) $n$-{\it ary Hom-superalgebra} is a triple ($A, \{\cdot , \cdots , \cdot \}, \alpha$) consisting of a 
superspace $A = A_0 \oplus A_1$, an $n$-linear map $\{\cdot , \cdots , \cdot \}: A^{\otimes n} \rightarrow A$ such that 
$\{A_i , \cdots , A_s \} \subset A_{i+ \cdots +s}$, and an even linear map $\alpha: A \rightarrow A$ such that
$\alpha (\{ x_1 , \cdots , x_n \}) = \{ \alpha (x_1), \cdots , \alpha (x_n) \}$ (multiplicativity).
\end{definition}

\par
One observes that if $\alpha = Id$ (the identity map), we get the corresponding definition of an $n$-{\it ary superalgebra}. One also notes
that $\overline{\alpha (x)} = \bar{x}$ for all homogeneous $x \in A$.
\par
We will be interested in binary ($n=2$), ternary ($n=3$) and binary-ternary Hom-superalgebras (i.e. Hom-superalgebras with binary and ternary
operations). For convenience, throughout this paper we assume that all Hom-(super)algebras are multiplicative.\\
\par
{\bf Definition 2.2.}
\cite{AM} Let ${\mathcal A} := (A, *, \alpha)$ be a binary Hom-superalgebra. The {\it supercommutator Hom-superalgebra}
(or the {\it minus Hom-superalgebra}) of $\mathcal A$ is the Hom-superalgebra ${\mathcal A}^{-} := (A, [\cdot , \cdot], \alpha)$, 
where $[x,y] := {1 \over 2}(x * y - (-1)^{\bar{x} \bar{y}} y * x)$ for all homogeneous $ x,y \in A$. The product $[\cdot , \cdot]$ is called the 
{\it supercommutator bracket}. The Hom-superalgebra ${\mathcal A}^{+} := (A, \circ , \alpha)$ is usually called
the {\it plus Hom-superalgebra} of ${\mathcal A}$.\\

\par
Throughout this paper and where there is no danger of confusion, we will also denote by ``$[\cdot , \cdot]$`` the binary operation in (Hom-)superalgebras. Besides the supercommutator of elements
in $\mathcal{A}$, one also considers the {\it super-Jordan product}
$$ x \circ y := \frac{1}{2}(xy + (-1)^{\bar{x} \bar{y}} yx),$$
the {\it Hom-Jordan associator}
$$as_{\mathcal{A}^{+}}(x,y,z) := (x \circ y) \circ \alpha (z) - \alpha (x) \circ (y \circ z) $$
and the {\it Hom-associator} 
$$as_{\mathcal{A}}(x,y,z) := (x*y) * \alpha (z) - \alpha (x) * (y*z) $$
for homogeneous $x,y,z$ in the given Hom-superalgebra \cite{MS}. When $\alpha = Id$ we recover the usual {\it Jordan associator} 
$as_{\mathcal{A}^{+}}(x,y,z) := (x \circ y) \circ z - x \circ (y \circ z)$ and 
{\it associator} $as_{\mathcal{A}}(x,y,z) := (x*y) * z - x * (y*z)$ respectively in usual superalgebras.
\par
The usual composition of maps  will also be denoted either by ''$\circ$`` or simply by juxtaposition.\\
\par
{\bf Definition 2.3.} 
A Hom-superalgebra ${\mathcal A} := (A, *, \alpha)$ is said to be {\it right Hom-alternative} if \\
\\
(2.1) \; $as_{\mathcal{A}}(x,y,z) = - (-1)^{\bar{y} \bar{z}} as_{\mathcal{A}}(x,z,y)$ ({\it right superalternativity}) \\
\\
for all homogeneous $x,y,z \in A$. \\

\par
Likewise is defined a left Hom-alternative superalgebra. A Hom- superalgebra that is both right and left Hom-alternative is said to be 
{\it Hom-alternative} \cite{AAM}.
\par
From (2.1), expanding Hom-associators, it is easily seen that (2.1) is equivalent to \\
\\
(2.2) \; $\alpha (x)(yz + (-1)^{\bar{y} \bar{z}} zy) = (xy)\alpha (z) + (-1)^{\bar{y} \bar{z}} (xz)\alpha (y)$. \\
\par
For $\alpha = Id$ in (2.1) we get a {\it right alternative superalgebra} $(A, *)$ that is first defined in \cite{ZS}
and (2.2) now reads as \\
\\
(2.3) \; $x(yz + (-1)^{\bar{y} \bar{z}} zy) = (xy)z + (-1)^{\bar{y} \bar{z}} (xz)y$. \\
\par
The following characterizations of right (Hom-)alternative superalgebras will be useful. \\
\par
{\bf Lemma 2.1.} 
{\it In any right alternative superalgebra the identity}\\
\\
(2.4) \; $as ([w,x],y,z) = [w, as(x,y,z)] + (-1)^{\bar{x}(\bar{y} + \bar{z})} [as(w,y,z),x]$ 
\par
\hspace{3.7cm}$- as(w,x, [y,z]) +  (-1)^{\bar{w}\bar{x}} as(x,w, [y,z])$ \\
\\
{\it holds for all homogeneous $w,x,y,z$}.\\

\par
{\it Proof.} The following identity in right alternative algebras is well-known ( see \cite{Thed}, identity (5)):
\\
(2.5) \; $as ([w,x],y,z) = [w, as(x,y,z)] + [as(w,y,z),x] - as(w,x, [y,z])$
\par
\hspace{3.8cm}$ + as(x,w, [y,z])$ \\
\\
(the reader is cautioned that, in the identity (5) of \cite{Thed}, there is a misprint on the third term where the
variables $z$ and $w$ are wrongly permuted, so it should be $[(x,w,z),y]$ instead of $[(x,z,w),y]$ as written). 
Applying the Kaplansky rule, the superization of (2.5) gives the identity (2.4). \hfill $\square$ \\ 
\par
The Hom-version of the identity (2.4) is given by the following \\
\par

{\bf Lemma 2.2.} {\it In any right Hom-alternative superalgebra the identity}\\
\\
(2.6) \; $as ([w,x], \alpha (y), \alpha (z)) =$ 
\par
\hspace{1.0cm}$[{\alpha}^{2} (w), as(x,y,z)] + (-1)^{\bar{x}(\bar{y} + \bar{z})} 
[as(w,y,z),{\alpha}^{2} (x)]$
\par
\hspace{1.0cm}$-as(\alpha (w),\alpha (x),[y,z])+(-1)^{\bar{w}\bar{x}} as(\alpha (x),\alpha (w),[y,z])$ \\
\\
{\it holds for all homogeneous $w,x,y,z$}.\\

\par
{\it Proof.} In \cite[Theorem 3.2]{Iss2} the following identity is proved to hold in any right Hom-alternative
superalgebra: \\
\\
(2.7) \; $as(wx, \alpha (y), \alpha (z)) = (-1)^{\bar{x}(\bar{y} + \bar{z})} as(w,y,z){\alpha}^{2}(x) + {\alpha}^{2}(w) as(x,y,z)$
\par
\hspace{4.2cm}$- 2 as(\alpha (w),\alpha (x),[y,z]) $ \\
\\
(the coefficient $2$ comes from the specific definition of the supercommutator). Let permute $w$ with $x$ in (2.7) to 
get \\
\\
(2.8) \; $(-1)^{\bar{w}\bar{x}} as(xw, \alpha (y), \alpha (z)) = (-1)^{\bar{w}(\bar{x} + \bar{y} 
+ \bar{z})} as(x,y,z){\alpha}^{2}(w) $
\par
\hspace{5.5cm}$+ (-1)^{\bar{w}\bar{x}} {\alpha}^{2}(x) as(w,y,z)$
\par
\hspace{5.5cm}$- 2 (-1)^{\bar{w}\bar{x}} as(\alpha (x),\alpha (w),[y,z])$. \\
\\
The subtraction of (2.8) from (2.7) leads to (2.6). \hfill $\square$  

\section{Right alternative superalgebras and Bol superalgebras}

In this section we prove that the supercommutator superalgebra of a right alternative superalgebra is a Bol superalgebra. First we recall the following\\
\par
{\bf Definition 3.1.} \cite{Ruk} A {\it Bol superalgebra} is a triple $(A, [\cdot , \cdot], \{ \cdot, \cdot, \cdot \})$ in which $A$ is a
superspace, $[\cdot , \cdot]$ and $\{ \cdot, \cdot, \cdot \}$ are binary and ternary operations respectively on $A$ 
such that \\
\par
(SB01) $[A_{i},A_{j}] \subseteq A_{i+j}$,
\par
(SB02) $\{ A_{i},A_{j},A_{k} \} \subseteq A_{i+j+k}$,
\par
(SB1) $[x,y] = - (-1)^{\bar{x} \bar{y}} [y,x]$,
\par
(SB2) $\{x,y,z\} = - (-1)^{\bar{x} \bar{y}} \{y,x,z\}$,
\par
(SB3) $\{x,y,z\} + (-1)^{\bar{x} (\bar{y} + \bar{z})} \{y,z,x\} + (-1)^{\bar{z} (\bar{x} + \bar{y})}\{z,x,y\} = 0$,
\par
(SB4) $\{x,y,[u,v]\} = [\{x,y,u\},v] + (-1)^{\bar{u} (\bar{x} + \bar{y})}[u,\{x,y,v\}]$
\par
\hspace{3.5cm}$+ (-1)^{(\bar{x} + \bar{y}) (\bar{u} + \bar{v})} (\{u,v,[x,y]\} - [[u,v],[x,y]])$,
\par
(SB5) $\{x,y,\{u,v,w\} \} = \{ \{x,y,u\}, v,w\} + (-1)^{\bar{u} (\bar{x} + \bar{y})} \{u, \{x,y,v\}, w\}$
\par
\hspace{4.2cm}$+ (-1)^{(\bar{x} + \bar{y}) (\bar{u} + \bar{v})} \{u,v, \{x,y,w\} \}$ \\
\\
for all homogeneous $u,v,w,x,y,z \in A$.\\
\par
Clearly, any Bol superalgebra with zero odd part is a (left) Bol algebra. If $[x,y] = 0$ for all homogeneous $x,y \in A$, then 
$(A, [\cdot , \cdot], \{ \cdot, \cdot, \cdot \})$ reduces to a {\it Lie supertriple system} $(A,\{ \cdot, \cdot, \cdot \})$
(see \cite{Tilg} for the definition of a {\it G-graded Lie triple system}).\\
\par
{\bf Example 3.1.} Let $A = A_0 \oplus A_1$ be a superspace where $A_0$ is a $2$-dimensional vector space with basis $\{i,j\}$ and  $A_1$ a 
$1$-dimensional vector space with basis $\{k\}$. Define on $A$ the following binary and ternary nonzero products:
\par
$[i,j]=j$, $[i,k]=k$,
\par
$[j,i]=-j$,
\par
$[k,i]=-k$, $[k,k]=j$,
\par
$\{i,j,i\}=-j$, $\{i,k,i\}=-k$,
\par
$\{j,i,i\}=j$,
\par
$\{k,i,i\}=k$. \\
Then it could be checked that $(A, [\cdot , \cdot], \{ \cdot, \cdot, \cdot \})$, with the multiplication table as above, is a ($3$-dimensional) 
Bol superalgebra. Observe that $(A, [\cdot , \cdot])$ is a $3$-dimensional Maltsev superalgebra \cite{AE} and the table for the ternary product 
as above is obtained using the $\mathbb{Z}_2$-graded version of the ternary product that produces a Bol algebra from a Maltsev 
algebra \cite{Mik1}.\\
\par
{\bf Definition 3.2.} \cite{Kac2} A {\it Jordan superalgebra} is a super vector space $V = V_0 \oplus V_1$ equipped
 with a bilinear product such that 
\par
(i) $V_{i} V_{j} \subset V_{i+j}$, $i,j \in \mathbb{Z}_2$,
\par
(ii) $xy = (-1)^{\bar{x}\bar{y}} yx$ ({\it supercommutativity}),
\par
(iii) ${\circlearrowleft}_{x,y,z} (-1)^{\bar{z} (\bar{x} + \bar{w})} as_{V}(xy,w,z)$ ({\it Jordan superidentity})\\
for all homogeneous $w,x,y,z \in V$, where ${\circlearrowleft}_{x,y,z}$ denotes the sum over cyclic permutation of
$x,y,z$. \\
\par
Consider in a right alternative superalgebra $(A,*)$ a ternary operation defined by \\
\\
(3.1) \; $[x,y,z] := 2 (-1)^{\bar{x} (\bar{y} + \bar{z})} as_{\mathcal{A}^{+}}(y,z,x)$. \\
\par
Another useful expression of the operation (3.1) is given by the following \\
\par

{\bf Lemma 3.1.} {\it If $(A,*)$ is a right alternative superalgebra, then} \\
\\
(3.2) \; $[x,y,z] := 2[[x,y],z] - (-1)^{\bar{z} (\bar{x} + \bar{y})} as_{\mathcal A}(z,x,y)$  \\
\\
{\it for all homogeneous $x,y,z \in A$}. \\

\par
{\it Proof.} We have 
\par
$as_{\mathcal{A}^{+}}(y,z,x) = (y \circ z) \circ x - y \circ (z \circ x)$
\par
\hspace{0.5cm}$= \frac{1}{2} \{ ( \frac{1}{2} (y*z + (-1)^{\bar{y}\bar{z}} z*y)*x
+ (-1)^{\bar{x} (\bar{y} + \bar{z})} x*( \frac{1}{2} (y*z + (-1)^{\bar{y}\bar{z}} z*y)) \}$
\par
\hspace{0.5cm}$- \frac{1}{2} \{ y* ( \frac{1}{2} (z*x + (-1)^{\bar{x}\bar{z}} x*z))
+ (-1)^{\bar{y} (\bar{x} + \bar{z})} ( \frac{1}{2} (z*x + (-1)^{\bar{x}\bar{z}} x*z)*y) \}$
\par
\hspace{0.5cm}$=\frac{1}{4} \{ (y*z + (-1)^{\bar{y}\bar{z}} z*y)*x 
+ (-1)^{\bar{x} (\bar{y} + \bar{z})} ((x*y)*z + (-1)^{\bar{y}\bar{z}} (x*z)*y)$
\par
\hspace{0.5cm}$- ((y*z)*x + (-1)^{\bar{x}\bar{z}} (y*x)*z) - (-1)^{\bar{y} (\bar{x} + \bar{z})}
(z*x + (-1)^{\bar{x}\bar{z}} x*z)*y \}$ (by (2.3))
\par
\hspace{0.5cm}$= \frac{1}{4} ( (-1)^{\bar{y}\bar{z}}(z*y)*x - (-1)^{\bar{y} (\bar{x} + \bar{z})} (z*x)*y 
+2 (-1)^{\bar{x} (\bar{y} + \bar{z})} [x,y]*z)$
\par
\hspace{0.5cm}$= \frac{1}{4} ((-1)^{\bar{y}\bar{z}} (z*y)*x - (-1)^{\bar{y} (\bar{x} + \bar{z})} (z*x)*y
+4 (-1)^{\bar{x} (\bar{y} + \bar{z})} [[x,y],z]$
\par
\hspace{0.5cm}$+ 2 (-1)^{\bar{y} (\bar{x} + \bar{z})} z*[x,y])$
\par
\hspace{0.5cm}$= \frac{1}{4} (4 (-1)^{\bar{x} (\bar{y} + \bar{z})} [[x,y],z] + (-1)^{\bar{y}\bar{z}} as_{\mathcal A}(z,y,x)
- (-1)^{\bar{y} (\bar{x} + \bar{z})} as_{\mathcal A}(z,x,y))$
\par
\hspace{0.5cm}$= \frac{1}{2} (2 (-1)^{\bar{x} (\bar{y} + \bar{z})} [[x,y],z] - (-1)^{\bar{y} (\bar{x} + \bar{z})} 
as_{\mathcal A}(z,x,y))$ (by the right superalternativity)\\
and then (3.2) follows. \hfill $\square$  
\\
\par
The next result is just the superization of the construction of a Lie triple system from any Jordan algebra \cite{Jac}. \\
\par

{\bf Proposition 3.1.} {\it Any Jordan superalgebra $(J, \circ )$ is a Lie supertriple system with respect to the 
operation $[x,y,z] := 2(x \circ (y \circ z) - (-1)^{\bar{x}\bar{y}} y \circ (x \circ z))$ on $(J, \circ )$}.\\

\par
{\it Proof.} The verification of identities (SB2) and (SB3) is straightforward. As for (SB5) we proceed as follows.
\par
First recall that the following identity holds in any Jordan superalgebra $(J, \circ )$:\\
\\
(3.3) \; $[[L(a),L(b)], L(c)] = L(a \circ (b \circ c)) - (-1)^{\bar{a}\bar{b}} L(b \circ (a \circ c))$, \\
\\
where $L(u)(v) := u \circ v$ (\cite{Kac2}). Moreover, $L(u) \in End(J)$, $\forall u \in J$ and $End(J)$ is a Lie 
superalgebra with respect to the superbracket $[f,g] = fg - (-1)^{\bar{f}\bar{g}} gf$ (the composition of two maps $f,g$
is written as $fg$). 
\par 
Now observe that $[x,y,z] = 2[L(x),L(y)](z)$ and so, if define the operator $D(x,y) := 2[L(x),L(y)]$, then
\par
$[x,y,[u,v,w]] - (-1)^{(\bar{x}+\bar{y}) (\bar{u}+\bar{v})} [u,v,[x,y,w]] = [D(x,y),D(u,v)](w)$ 
\par 
$= 4[[L(x),L(y)],[L(u),L(v)]](w)$. \\
Next we compute 
\par
$[[x,y,u],v,w] = 2[L_{D(x,y)(u)}, L(v)](w)$
\par
\hspace{2.5cm}$= 4[L(x(yu)) - (-1)^{\bar{x}\bar{y}} L(y(xu)), L(v)](w)$
\par
\hspace{2.5cm}$= 4 [[[L(x),L(y)],L(u)],L(v)](w)$ (by (3.3)) \\
and, likewise, 
\par
$(-1)^{\bar{u} (\bar{x} + \bar{y})} [u,[x,y,v],w] = -4 (-1)^{\bar{u}\bar{v}} [[[L(x),L(y)],L(v)],L(u)](w)$. \\
Thus 
\par
$[[x,y,u],v,w] + (-1)^{\bar{u} (\bar{x} + \bar{y})} [u,[x,y,v],w]$
\par
$= 4( [[[L(x),L(y)],L(u)],L(v)] - (-1)^{\bar{u}\bar{v}} [[[L(x),L(y)],L(v)],L(u)]) (w)$
\par
$= 4( [[[L(x),L(y)],L(u)],L(v)] + (-1)^{\bar{v}(\bar{u} + \bar{x} + \bar{y})} [[L(v),[L(x),L(y)]],L(u)])(w)$
\par
$= -4 (-1)^{(\bar{x}+\bar{y}) (\bar{u}+\bar{v})} [[L(u),L(v)],[L(x),L(y)]](w)$ (by the Jacobi superidentity in
$End(J)$)
\par
$= [2[L(x),L(y)],2[L(u),L(v)]](w)$
\par
$= [D(x,y),D(u,v)](w)$
\par
$= [x,y,[u,v,w]] - (-1)^{(\bar{x}+\bar{y}) (\bar{u}+\bar{v})}[u,v,[x,y,w]]$ \\
which proves (SB5). This completes the proof. \hfill $\square$  
\\
\par
In \cite{Mik2} it is proved that on any right alternative algebra one may define a Bol algebra structure. The $\mathbb{Z}_{2}$-graded version of this 
result is given by the following \\
\par

{\bf Theorem 3.1.} {\it The supercommutator algebra of any right alternative superalgebra is a Bol superalgebra}. \\

\par
{\it Proof.} Let $A$ be a right alternative superalgebra. Then $(A, \circ )$ is a Jordan superalgebra \cite{Sh2}
and Proposition 3.1 implies that $(A, [\cdot, \cdot, \cdot])$ is a Lie supertriple system, where ''$[\cdot, \cdot, \cdot]$``
is defined by (3.1) (or (3.2)). Therefore $(A, \{ \cdot, \cdot, \cdot \})$, with $\{x,y,z\} := \frac{1}{2} [x,y,z]$, 
is also a Lie supertriple system. Next we are done if we prove that (SB4) holds since (SB1) follows from the 
definition of the supercommutator.
\par
From (3.2) we have\\
\\
(3.4) \; $\frac{1}{2} as_{\mathcal A}(z,x,y) = (-1)^{\bar{z} (\bar{x} + \bar{y})} ([[x,y],z] - \{x,y,z\})$. \\
\\
Then, using (3.4), the identity (2.4) implies that
\par
$(-1)^{(\bar{x}+\bar{y}) (\bar{u}+\bar{v})} ([[x,y],[u,v]]-\{x,y,[u,v]\})$
\par
$= [u, (-1)^{\bar{v} (\bar{x} + \bar{y})} ([[x,y],v] - \{x,y,v\})]$ 
\par
$+ (-1)^{\bar{v} (\bar{x} + \bar{y})} [(-1)^{\bar{u} (\bar{x} + \bar{y})} ([[x,y],u] - \{x,y,u\}), v]$
\par
$- (-1)^{\bar{u} (\bar{v}+\bar{x} +\bar{y})} ([[v,[x,y]],u] - \{v,[x,y],u\})$
\par
$+ (-1)^{\bar{u}\bar{v}} ((-1)^{\bar{v} (\bar{u}+\bar{x} +\bar{y})} ([[u,[x,y]],v] - \{u,[x,y],v\}))$ \\
i.e. 
\par
$(-1)^{(\bar{x}+\bar{y}) (\bar{u}+\bar{v})} \{x,y,[u,v]\}$ 
\par
$= (-1)^{(\bar{x}+\bar{y}) (\bar{u}+\bar{v})} [\{x,y,u\},v]
+ (-1)^{\bar{v} (\bar{x} + \bar{y})} [u, \{x,y,v\}]$
\par
$- (-1)^{(\bar{x}+\bar{y}) (\bar{u}+\bar{v})} \{[x,y],u,v\} + (-1)^{\bar{u}\bar{v} +(\bar{x}+\bar{y}) (\bar{u}+\bar{v})}
\{[x,y],v,u\}$
\par
$+ (-1)^{(\bar{x}+\bar{y}) (\bar{u}+\bar{v})} [[x,y],[u,v]]$\\
or\\
\\
(3.5) \; $\{x,y,[u,v]\} = [\{x,y,u\},v] + (-1)^{\bar{u} (\bar{x} + \bar{y})} [u, \{x,y,v\}]$
\par
\hspace{3.0cm}$ - \{[x,y],u,v\} + (-1)^{\bar{u}\bar{v}} \{[x,y],v,u\} + [[x,y],[u,v]]$. \\
\\
Now in (3.5) observe that the expression $- \{[x,y],u,v\} + (-1)^{\bar{u}\bar{v}} \{[x,y],v,u\}$ can be transformed as
follows.
\par
$- \{[x,y],u,v\} + (-1)^{\bar{u}\bar{v}} \{[x,y],v,u\} = (-1)^{\bar{u} (\bar{x} + \bar{y})} \{u,[x,y],v\}$
\par
$+ (-1)^{\bar{u}\bar{v}} \{[x,y],v,u\}$
\par
$= - (-1)^{\bar{u}\bar{v} +(\bar{x}+\bar{y}) (\bar{u}+\bar{v})} \{v,u,[x,y]\}$ (by (SB3)). \\
Therefore (3.5) now reads as
\par
$\{x,y,[u,v]\} = [\{x,y,u\},v] + (-1)^{\bar{u} (\bar{x} + \bar{y})} [u, \{x,y,v\}]$
\par
\hspace{2.5cm}$+ (-1)^{(\bar{x}+\bar{y}) (\bar{u}+\bar{v})} \{u,v,[x,y]\} + [[x,y],[u,v]]$ \\
and so we get (SB4). This completes the proof. \hfill $\square$  
\\ 
\par
An application of Theorem 3.1 is given in Example 5.1 below.

\section{Hom-Jordan-admissibility of right Hom-alternative superalgebras}

In \cite[Proposition 1]{Sh2} it is proved that every right alternative superalgebra is Jordan-admissible and the ungraded
Hom-version of this result is considered in \cite{Yau2} where it is specifically proved that every multiplicative right
Hom-alternative algebra is Hom-Jordan-admissible (\cite[Theorem 4.3]{Yau2}; see also \cite[Theorem 5.6]{Yau4}, for the 
case of Hom-alternative algebras). The $\mathbb{Z}_2$-graded extension of \cite[Theorem 5.6]{Yau4} to Hom-superalgebras
is considered in \cite{AAM}.
\par
In this section we consider the $\mathbb{Z}_2$-graded generalization of Theorem 4.3 of \cite{Yau2}. Since the twisting
map of a given Hom-algebra is neither injective nor surjective in general, the method used to prove Proposition 1 in
\cite{Sh2} cannot be reported to the Hom-algebra setting at the present stage of the theory of Hom-algebras. So we
proceed otherwise (see Lemma 4.1 and Theorem 4.1 below) by considering the $\mathbb{Z}_2$-graded generalization of the 
method used in \cite{Yau4} for a proof of the Hom-Jordan admissibility of Hom-alternative algebras. Observe that such
a generalization is first mentioned in \cite[Theorem 6.1 and Lemma 6.1]{AAM} but we point out that, in the formulation
of Lemma 6.1 of \cite{AAM} and its proof, some terms are repeated and others are missing. For completeness and because
of its key role in our setting, we give a full formulation and a proof of this lemma in the present paper (see Lemma 4.1)
and next prove that every multiplicative right Hom-alternative superalgebra is Hom-Jordan-admissible (Theorem 4.1).
\par
First we recall some relevant definitions from \cite{AAM}.\\
\par
{\bf Definition 4.1.} \cite{AAM} (i) A Hom-superalgebra $\mathcal{A}$ is called a {\it Hom-Jordan superalgebra} if it
is supercommutative (i.e. $xy = (-1)^{\bar{x} \bar{y}} yx$) and the {\it Hom-Jordan superidentity}
\par
${\circlearrowleft}_{x,y,t} \; (-1)^{\bar{t}(\bar{x} + \bar{z})} as_{\mathcal{A}} (xy, \alpha (z), \alpha (t)) = 0$\\
holds for all homogeneous $x,y,z,t$.
\par
(ii) A Hom-superalgebra $\mathcal{A} = (A,*,\alpha)$ is said to be Hom-Jordan-admissible if its plus Hom-superalgebra
${\mathcal{A}}^{+} = (A, \circ , \alpha)$ is a Hom-Jordan superalgebra. \\
\par
For $\alpha = Id$ we get the notions of a Jordan superalgebra and a Jordan-admissible superalgebra respectively. \\
\par
{\bf Remark 4.1.} Since the ''$\circ$`` is supercommutative, a Hom-superalgebra $\mathcal{A} = (A,*,\alpha)$ is 
Hom-Jordan-admissible if and only if ${\mathcal{A}}^{+}$ satisfies the Hom-Jordan superidentity\\
\\
(4.1) \; ${\circlearrowleft}_{x,y,t} \; (-1)^{\bar{t}(\bar{x} + \bar{z})} as_{\mathcal{A}^{+}} (x \circ y, \alpha (z), \alpha (t)) 
= 0$ \\
\\
for all homogeneous $x,y,z,t \in A$.\\
\par
Expanding the Hom-associator $as_{\mathcal{A}^{+}}$ in (4.1) and using the supercommutativity of ''$\circ$``, one checks 
that (4.1) is equivalent to \\
\\
(4.2) \; ${\circlearrowleft}_{x,y,t} \; (-1)^{\bar{x}(\bar{y} + \bar{z}) + \bar{y}(\bar{t} + \bar{z})} 
\{ (\alpha (t) \circ \alpha (z)) \circ \alpha (y \circ x) - {\alpha}^{2} (t) \circ ( \alpha (z) \circ (y \circ x)) \}
=0$. \\
\par
Below the operation ''$*$`` is denoted by juxtaposition. \\
\par

{\bf Lemma 4.1.} {\it Let $\mathcal{A} = (A,*,\alpha)$ be any Hom-superalgebra and ${\mathcal{A}}^{+} = (A, \circ , \alpha)$
its plus Hom-superalgebra. Then}\\
\\
(4.3) \; $8 {\circlearrowleft}_{x,y,t} \; (-1)^{\bar{x}(\bar{y} + \bar{z}) + \bar{y}(\bar{t} + \bar{z})} 
\{ (\alpha (t) \circ \alpha (z)) \circ \alpha (y \circ x) - {\alpha}^{2} (t) \circ ( \alpha (z) \circ (y \circ x)) \}$
\par
$= {\circlearrowleft}_{x,y,t} \; (-1)^{\bar{x}(\bar{y} + \bar{z}) + \bar{y}(\bar{t} + \bar{z})} 
\{ as_{\mathcal{A}} (\alpha (t),\alpha (z),yx) + (-1)^{\bar{x} \bar{y}} as_{\mathcal{A}} (\alpha (t),\alpha (z),xy)$
\par
$-(-1)^{\bar{z}(\bar{x} + \bar{y}) + \bar{t}(\bar{x} + \bar{y} + \bar{z})} as_{\mathcal{A}} (yx,\alpha (z),\alpha (t))$
\par
$-(-1)^{\bar{x}\bar{y} + \bar{z}(\bar{x} + \bar{y}) + \bar{t}(\bar{x} + \bar{y} + \bar{z})} as_{\mathcal{A}} (xy,\alpha (z),\alpha (t))$
\par
$+  (-1)^{\bar{t} \bar{z}} as_{\mathcal{A}} (\alpha (z),\alpha (t),yx) + (-1)^{\bar{t}\bar{z}+\bar{x}\bar{y}}
as_{\mathcal{A}} (\alpha (z),\alpha (t),xy)$
\par
$- (-1)^{(\bar{t}+\bar{z})(\bar{x}+\bar{y})} as_{\mathcal{A}} (yx,\alpha (t),\alpha (z)) 
- (-1)^{\bar{x}\bar{y}+(\bar{t}+\bar{z})(\bar{x}+\bar{y})} as_{\mathcal{A}} (xy,\alpha (t),\alpha (z))$
\par
$+ (-1)^{\bar{z}(\bar{x} + \bar{y})} as_{\mathcal{A}} (\alpha (t),yx,\alpha (z))
+ (-1)^{\bar{x}\bar{y} + \bar{z}(\bar{x} + \bar{y})} as_{\mathcal{A}} (\alpha (t),xy,\alpha (z))$
\par
$-(-1)^{\bar{t}(\bar{x} + \bar{y} + \bar{z})} as_{\mathcal{A}} (\alpha (z),yx,\alpha (t))
-(-1)^{\bar{x}\bar{y} +\bar{t}(\bar{x} + \bar{y} + \bar{z})} as_{\mathcal{A}} (\alpha (z),xy,\alpha (t)) \}$
\par
$- (-1)^{\bar{z}(\bar{t} + \bar{x}) + \bar{y}(\bar{t} + \bar{z})} [{\alpha}^{2} (z), as_{\mathcal{A}} (t,x,y)]
- (-1)^{\bar{z}(\bar{x} + \bar{y}) + \bar{t}(\bar{x} + \bar{z})} [{\alpha}^{2} (z), as_{\mathcal{A}} (x,y,t)]$
\par
$- (-1)^{\bar{z}(\bar{y} + \bar{t}) + \bar{x}(\bar{y} + \bar{z})} [{\alpha}^{2} (z), as_{\mathcal{A}} (y,t,x)]
- (-1)^{\bar{z}(\bar{t}+\bar{x}+\bar{y})+ \bar{y}(\bar{t}+\bar{x})} [{\alpha}^{2} (z), as_{\mathcal{A}} (t,y,x)]$
\par
$- (-1)^{\bar{z}(\bar{t}+\bar{x}+\bar{y})+ \bar{x}(\bar{t}+\bar{y})} [{\alpha}^{2} (z), as_{\mathcal{A}} (y,x,t)]
-(-1)^{\bar{z}(\bar{t}+\bar{x}+\bar{y})+\bar{t}(\bar{x}+\bar{y})}[{\alpha}^{2}(z),as_{\mathcal{A}}(x,t,y)]$ \\
\\
{\it for all homogeneous $t,x,y,z \in A$}. \\

\par
{\it Proof.} Observe that 
\par
$(\alpha (t) \circ \alpha (z)) \circ \alpha (y \circ x)  = \frac{1}{2} \{ (\alpha (t) \circ \alpha (z))\alpha (y \circ x)$
\par
$+ (-1)^{(\bar{t}+\bar{z})(\bar{x}+\bar{y})} \alpha (y \circ x) (\alpha (t) \circ \alpha (z)) \}$
\par
$= \frac{1}{8} \{ (\alpha (t)\alpha (z))\alpha (yx) + (-1)^{\bar{x}\bar{y}}(\alpha (t)\alpha (z))\alpha (xy) 
+ (-1)^{\bar{t}\bar{z}}(\alpha (z)\alpha (t))\alpha (yx)$
\par
$ +(-1)^{\bar{t}\bar{z} +\bar{x}\bar{y}}(\alpha (z)\alpha (t))\alpha (xy) + (-1)^{(\bar{t}+\bar{z})(\bar{x}+\bar{y})} 
\alpha (yx)(\alpha (t)\alpha (z))$
\par
$ + (-1)^{\bar{t}\bar{z}+(\bar{t}+\bar{z})(\bar{x}+\bar{y})}\alpha (yx)(\alpha (z)\alpha (t))
+ (-1)^{\bar{x}\bar{y}+(\bar{t}+\bar{z})(\bar{x}+\bar{y})}\alpha (xy)(\alpha (t)\alpha (z))$
\par
$(-1)^{\bar{t}\bar{z} +\bar{x}\bar{y}+(\bar{t}+\bar{z})(\bar{x}+\bar{y})}
\alpha (xy) (\alpha (z)\alpha (t)) \}$. \\
Likewise we have
\par
${\alpha}^{2} (t) \circ ( \alpha (z) \circ (y \circ x)) = \frac{1}{8} \{ {\alpha}^{2}(t) (\alpha (z)(yx)) 
+ (-1)^{\bar{x}\bar{y}}{\alpha}^{2}(t) (\alpha (z)(xy))$
\par
$+ (-1)^{\bar{z}(\bar{x} + \bar{y})}{\alpha}^{2}(t) ((yx)\alpha (z))
+ (-1)^{\bar{x}\bar{y} +\bar{z}(\bar{x} + \bar{y})}{\alpha}^{2}(t) ((xy)\alpha (z))$
\par
$+ (-1)^{\bar{t}(\bar{x} + \bar{y} + \bar{z})}(\alpha (z)(yx)) {\alpha}^{2}(t) 
+ (-1)^{\bar{x}\bar{y} + \bar{t}(\bar{x} + \bar{y} + \bar{z})}(\alpha (z)(xy)) {\alpha}^{2}(t)$
\par
$+ (-1)^{\bar{z}(\bar{x} + \bar{y}) +\bar{t}(\bar{x} + \bar{y} + \bar{z})}((yx)\alpha (z)) {\alpha}^{2}(t)
+(-1)^{\bar{x}\bar{y} +\bar{z}(\bar{x} + \bar{y}) +\bar{t}(\bar{x}+\bar{y}+\bar{z})}((xy)\alpha (z)){\alpha}^{2}(t) \}$.\\
Therefore, replacing $(\alpha (t) \circ \alpha (z)) \circ \alpha (y \circ x)$ and
${\alpha}^{2} (t) \circ ( \alpha (z) \circ (y \circ x))$ with their corresponding expressions as above and next
rearranging terms, we get
\par
$8  {\circlearrowleft}_{x,y,t} \; (-1)^{\bar{x}(\bar{y} + \bar{z}) + \bar{y}(\bar{t} + \bar{z})} 
\{ (\alpha (t) \circ \alpha (z)) \circ \alpha (y \circ x) - {\alpha}^{2} (t) \circ ( \alpha (z) \circ (y \circ x)) \}$
\par
$= {\circlearrowleft}_{x,y,t} \; (-1)^{\bar{x}(\bar{y} + \bar{z}) + \bar{y}(\bar{t} + \bar{z})} 
\{ as_{\mathcal{A}} (\alpha (t), \alpha (z), yx) + (-1)^{\bar{x}\bar{y}} as_{\mathcal{A}}(\alpha(t),\alpha (z), xy)$
\par
$- (-1)^{\bar{z}(\bar{x}+\bar{y}) + \bar{t}(\bar{x}+\bar{y}+\bar{z})}as_{\mathcal{A}} (yx,\alpha(z),\alpha (t))$
\par
$- (-1)^{\bar{x}\bar{y}+\bar{z}(\bar{x}+\bar{y})+\bar{t}(\bar{x}+\bar{y}+\bar{z})}as_{\mathcal{A}} (xy,\alpha(z),\alpha (t))$
\par
$+ (-1)^{\bar{t}\bar{z}} (\alpha (z)\alpha (t))\alpha (yx) + (-1)^{\bar{t}\bar{z} +\bar{x}\bar{y}}(\alpha (z)\alpha (t))\alpha (xy)$
\par
$+ (-1)^{(\bar{t}+\bar{z})(\bar{x}+\bar{y})} \alpha (yx)(\alpha (t)\alpha (z))
+ (-1)^{\bar{x}\bar{y}+(\bar{t}+\bar{z})(\bar{x}+\bar{y})}\alpha (xy)(\alpha (t)\alpha (z))$
\par
$- (-1)^{\bar{z}(\bar{x} + \bar{y})}{\alpha}^{2}(t) ((yx)\alpha (z))
- (-1)^{\bar{x}\bar{y} +\bar{z}(\bar{x} + \bar{y})}{\alpha}^{2}(t) ((xy)\alpha (z))$
\par
$- (-1)^{\bar{t}(\bar{x} + \bar{y} + \bar{z})}(\alpha (z)(yx)) {\alpha}^{2}(t)
- (-1)^{\bar{x}\bar{y} + \bar{t}(\bar{x} + \bar{y} + \bar{z})}(\alpha (z)(xy)) {\alpha}^{2}(t) \}$. \\
From the definition of the Hom-associator $as_{\mathcal{A}}$, each of the last eight terms of the right-hand side of
the expression above is written as follows.
\par
$\bullet$ $(-1)^{\bar{t}\bar{z}} (\alpha (z)\alpha (t))\alpha (yx)= (-1)^{\bar{t}\bar{z}} as_{\mathcal{A}}(\alpha(z),\alpha (t),yx) 
+ (-1)^{\bar{t}\bar{z}}{\alpha}^{2}(z)( \alpha (t) (yx))$;
\par
$\bullet$ $(-1)^{\bar{t}\bar{z} +\bar{x}\bar{y}}(\alpha (z)\alpha (t))\alpha (xy)
= (-1)^{\bar{t}\bar{z} +\bar{x}\bar{y}}as_{\mathcal{A}}(\alpha(z),\alpha (t),xy)$
\par
$+(-1)^{\bar{t}\bar{z} +\bar{x}\bar{y}} {\alpha}^{2}(z)( \alpha (t) (xy))$;
\par
$\bullet$ $(-1)^{(\bar{t}+\bar{z})(\bar{x}+\bar{y})} \alpha (yx)(\alpha (t)\alpha (z))
= - (-1)^{(\bar{t}+\bar{z})(\bar{x}+\bar{y})} as_{\mathcal{A}} (yx, \alpha (t), \alpha (z))$
\par
\hspace{6cm}$+ (-1)^{(\bar{t}+\bar{z})(\bar{x}+\bar{y})} ((yx)\alpha (t)) {\alpha}^{2}(z)$;
\par
$\bullet$ $(-1)^{\bar{x}\bar{y}+(\bar{t}+\bar{z})(\bar{x}+\bar{y})}\alpha (xy)(\alpha (t)\alpha (z))
= -(-1)^{\bar{x}\bar{y}+(\bar{t}+\bar{z})(\bar{x}+\bar{y})} as_{\mathcal{A}} (xy, \alpha (t), \alpha (z))$
\par
\hspace{6.5cm}$+ (-1)^{\bar{x}\bar{y}+(\bar{t}+\bar{z})(\bar{x}+\bar{y})} ((xy)\alpha (t)) {\alpha}^{2}(z)$;
\par
$\bullet$ $- (-1)^{\bar{z}(\bar{x} + \bar{y})}{\alpha}^{2}(t) ((yx)\alpha (z))
= (-1)^{\bar{z}(\bar{x} + \bar{y})} as_{\mathcal{A}} (\alpha (t),yx,\alpha (z))$
\par
\hspace{5.5cm}$- (-1)^{\bar{z}(\bar{x} + \bar{y})} (\alpha (t) (yx)) {\alpha}^{2}(z)$;
\par
$\bullet$ $- (-1)^{\bar{x}\bar{y} + \bar{z}(\bar{x} + \bar{y})}{\alpha}^{2}(t) ((xy)\alpha (z))
= (-1)^{\bar{x}\bar{y} + \bar{z}(\bar{x} + \bar{y})} as_{\mathcal{A}} (\alpha (t),xy,\alpha (z))$
\par
\hspace{6cm}$- (-1)^{\bar{x}\bar{y} + \bar{z}(\bar{x} + \bar{y})} (\alpha (t) (xy)) {\alpha}^{2}(z)$;
\par
$\bullet$ $- (-1)^{\bar{t}(\bar{x} + \bar{y} + \bar{z})}(\alpha (z)(yx)) {\alpha}^{2}(t)
= - (-1)^{\bar{t}(\bar{x} + \bar{y} + \bar{z})}as_{\mathcal{A}} (\alpha (z),yx,\alpha (t))$
\par
\hspace{6cm}$ -(-1)^{\bar{t}(\bar{x} + \bar{y} + \bar{z})} {\alpha}^{2}(z) ((yx)\alpha (t))$;
\par
$\bullet$ $- (-1)^{\bar{x}\bar{y} + \bar{t}(\bar{x} + \bar{y} + \bar{z})}(\alpha (z)(xy)) {\alpha}^{2}(t)
= - (-1)^{\bar{x}\bar{y} + \bar{t}(\bar{x} + \bar{y} + \bar{z})}as_{\mathcal{A}} (\alpha (z),xy,\alpha (t))$
\par
\hspace{6.5cm}$ -(-1)^{\bar{x}\bar{y} +\bar{t}(\bar{x} + \bar{y} + \bar{z})} {\alpha}^{2}(z) ((xy)\alpha (t))$.
\par
So we get \\
\\
(4.4) \; $8  {\circlearrowleft}_{x,y,t} \; (-1)^{\bar{x}(\bar{y} + \bar{z}) + \bar{y}(\bar{t} + \bar{z})} 
\{ (\alpha (t) \circ \alpha (z)) \circ \alpha (y \circ x) - {\alpha}^{2} (t) \circ ( \alpha (z) \circ (y \circ x)) \}$
\par
$= {\circlearrowleft}_{x,y,t} \; (-1)^{\bar{x}(\bar{y} + \bar{z}) + \bar{y}(\bar{t} + \bar{z})} 
\{ as_{\mathcal{A}} (\alpha (t),\alpha (z),yx) + (-1)^{\bar{x}\bar{y}} as_{\mathcal{A}} (\alpha (t),\alpha (z),xy)$
\par
$-(-1)^{\bar{z}(\bar{x} + \bar{y}) + \bar{t}(\bar{x} + \bar{y} + \bar{z})}as_{\mathcal{A}}(yx,\alpha (z),\alpha (t))$
\par
$-(-1)^{\bar{x}\bar{y}+\bar{z}(\bar{x}+\bar{y})+\bar{t}(\bar{x}+\bar{y}+\bar{z})}as_{\mathcal{A}}(xy,\alpha (z),\alpha (t))$
\par
$+(-1)^{\bar{t}\bar{z}} as_{\mathcal{A}}(\alpha(z),\alpha (t),yx) 
+ (-1)^{\bar{t}\bar{z} +\bar{x}\bar{y}}as_{\mathcal{A}}(\alpha(z),\alpha (t),xy)$
\par
$- (-1)^{(\bar{t}+\bar{z})(\bar{x}+\bar{y})} as_{\mathcal{A}} (yx, \alpha (t), \alpha (z))
- (-1)^{\bar{x}\bar{y}+(\bar{t}+\bar{z})(\bar{x}+\bar{y})} as_{\mathcal{A}} (xy, \alpha (t), \alpha (z))$
\par
$+ (-1)^{\bar{z}(\bar{x} + \bar{y})} as_{\mathcal{A}} (\alpha (t),yx,\alpha (z)) 
+ (-1)^{\bar{x}\bar{y}+\bar{z}(\bar{x}+\bar{y})} as_{\mathcal{A}} (\alpha (t),xy,\alpha (z))$
\par
$- (-1)^{\bar{t}(\bar{x} + \bar{y} + \bar{z})} as_{\mathcal{A}} (\alpha (z),yx,\alpha (t))
- (-1)^{\bar{x}\bar{y} + \bar{t}(\bar{x} + \bar{y} + \bar{z})} as_{\mathcal{A}} (\alpha (z),xy,\alpha (t))$
\par
$+(-1)^{\bar{t}\bar{z}}{\alpha}^{2}(z)(\alpha (t)(yx))+(-1)^{\bar{t}\bar{z}+\bar{x}\bar{y}}{\alpha}^{2}(z)(\alpha (t)(xy))$
\par
$+ (-1)^{(\bar{t}+\bar{z})(\bar{x}+\bar{y})} ((yx)\alpha (t)) {\alpha}^{2}(z)
+ (-1)^{\bar{x}\bar{y} + (\bar{t}+\bar{z})(\bar{x}+\bar{y})} ((xy)\alpha (t)) {\alpha}^{2}(z)$
\par
$- (-1)^{\bar{z}(\bar{x} + \bar{y})} (\alpha (t) (yx)) {\alpha}^{2}(z)
- (-1)^{\bar{x}\bar{y} + \bar{z}(\bar{x} + \bar{y})} (\alpha (t) (xy)) {\alpha}^{2}(z)$
\par
$- (-1)^{\bar{t}(\bar{x} + \bar{y} + \bar{z})} {\alpha}^{2}(z) ((yx)\alpha (t))
- (-1)^{\bar{x}\bar{y} + \bar{t}(\bar{x} + \bar{y} + \bar{z})} {\alpha}^{2}(z) ((xy)\alpha (t)) \}$. \\
\\
Next we transform the sum of the last eight terms in the right-hand side of (4.4) as follows. 
\par
${\circlearrowleft}_{x,y,t} \; (-1)^{\bar{x}(\bar{y} + \bar{z}) + \bar{y}(\bar{t} + \bar{z})} 
\{ (-1)^{\bar{t}\bar{z}}{\alpha}^{2}(z)(\alpha (t)(yx))+(-1)^{\bar{t}\bar{z}+\bar{x}\bar{y}}{\alpha}^{2}(z)(\alpha (t)(xy))$
\par
$+ (-1)^{(\bar{t}+\bar{z})(\bar{x}+\bar{y})} ((yx)\alpha (t)) {\alpha}^{2}(z)
+ (-1)^{\bar{x}\bar{y} + (\bar{t}+\bar{z})(\bar{x}+\bar{y})} ((xy)\alpha (t)) {\alpha}^{2}(z)$
\par
$- (-1)^{\bar{z}(\bar{x} + \bar{y})} (\alpha (t) (yx)) {\alpha}^{2}(z)
- (-1)^{\bar{x}\bar{y} + \bar{z}(\bar{x} + \bar{y})} (\alpha (t) (xy)) {\alpha}^{2}(z)$
\par
$- (-1)^{\bar{t}(\bar{x} + \bar{y} + \bar{z})} {\alpha}^{2}(z) ((yx)\alpha (t))
- (-1)^{\bar{x}\bar{y} + \bar{t}(\bar{x} + \bar{y} + \bar{z})} {\alpha}^{2}(z) ((xy)\alpha (t)) \}$
\par
$= (-1)^{\bar{x}(\bar{y} + \bar{z}) + \bar{y}(\bar{t} + \bar{z}) + \bar{t}\bar{z}} {\alpha}^{2}(z)(\alpha (t)(yx))
+ (-1)^{\bar{t}(\bar{x} + \bar{z}) + \bar{x}(\bar{y} + \bar{z}) + \bar{y}\bar{z}} {\alpha}^{2}(z)(\alpha (y)(xt))$
\par
$+ (-1)^{\bar{y}(\bar{t} + \bar{z}) + \bar{t}(\bar{x} + \bar{z}) + \bar{x}\bar{z}} {\alpha}^{2}(z)(\alpha (x)(ty))
+ (-1)^{\bar{z}(\bar{t} + \bar{x}) + \bar{y}(\bar{t} + \bar{z})} {\alpha}^{2}(z)(\alpha (t)(xy))$
\par
$+ (-1)^{\bar{z}(\bar{x} + \bar{y}) + \bar{t}(\bar{x} + \bar{z})} {\alpha}^{2}(z)(\alpha (x)(yt))
+ (-1)^{\bar{z}(\bar{y} + \bar{t}) + \bar{x}(\bar{y} + \bar{z})} {\alpha}^{2}(z)(\alpha (y)(tx))$
\par
$+ (-1)^{\bar{x}(\bar{t} + \bar{y})} ((yx) \alpha (t)) {\alpha}^{2}(z)
+ (-1)^{\bar{t}(\bar{y} + \bar{x})} ((xt) \alpha (y)) {\alpha}^{2}(z)$
\par
$+ (-1)^{\bar{y}(\bar{x} + \bar{t})} ((ty) \alpha (x)) {\alpha}^{2}(z)$
\par
$+ (-1)^{\bar{t}\bar{x}} ((xy) \alpha (t)) {\alpha}^{2}(z)
+ (-1)^{\bar{x}\bar{y}} ((yt) \alpha (x)) {\alpha}^{2}(z)
+ (-1)^{\bar{y}\bar{t}} ((tx) \alpha (y)) {\alpha}^{2}(z)$
\par
$- (-1)^{\bar{y}(\bar{t} + \bar{x})} (\alpha (t) (yx)) {\alpha}^{2}(z)
- (-1)^{\bar{x}(\bar{y} + \bar{t})} (\alpha (y) (xt)) {\alpha}^{2}(z)$
\par
$- (-1)^{\bar{t}(\bar{x} + \bar{y})} (\alpha (x) (ty)) {\alpha}^{2}(z)$
\par
$- (-1)^{\bar{y}\bar{t}} (\alpha (t) (xy)) {\alpha}^{2}(z)
- (-1)^{\bar{t}\bar{x}} (\alpha (x) (yt)) {\alpha}^{2}(z)
- (-1)^{\bar{x}\bar{y}} (\alpha (y) (tx)) {\alpha}^{2}(z)$
\par
$- (-1)^{\bar{x}(\bar{y} + \bar{z}) + \bar{t}(\bar{x} + \bar{z}) + \bar{y}\bar{z}} {\alpha}^{2}(z) ((yx) \alpha (t))
- (-1)^{\bar{t}(\bar{x} + \bar{z}) + \bar{y}(\bar{t} + \bar{z}) + \bar{x}\bar{z}} {\alpha}^{2}(z) ((xt) \alpha (y))$
\par
$-(-1)^{\bar{y}(\bar{t} + \bar{z}) + \bar{x}(\bar{y} + \bar{z}) + \bar{t}\bar{z}} {\alpha}^{2}(z) ((ty) \alpha (x))
- (-1)^{\bar{x}(\bar{t} + \bar{z}) + \bar{z}(\bar{t} + \bar{y})} {\alpha}^{2}(z) ((xy) \alpha (t))$
\par
$- (-1)^{\bar{y}(\bar{x} + \bar{z}) + \bar{z}(\bar{x} + \bar{t})} {\alpha}^{2}(z) ((yt) \alpha (x))
- (-1)^{\bar{t}(\bar{y} + \bar{z}) + \bar{z}(\bar{y} + \bar{x})} {\alpha}^{2}(z) ((tx) \alpha (y))$
\par
$= - (-1)^{\bar{z}(\bar{t} + \bar{x} + \bar{y}) + \bar{y}(\bar{t} + \bar{x})} [{\alpha}^{2}(z), as(t,y,x)]
+ (-1)^{\bar{z}(\bar{t} + \bar{x} + \bar{y}) + \bar{y}(\bar{t} + \bar{x})} [{\alpha}^{2}(z), (ty) \alpha (x)]$
\par
$ - (-1)^{\bar{z}(\bar{t} + \bar{x} + \bar{y}) + \bar{x}(\bar{t} + \bar{y})} [{\alpha}^{2}(z), as(y,x,t)]
+ (-1)^{\bar{z}(\bar{t} + \bar{x} + \bar{y}) + \bar{x}(\bar{t} + \bar{y})} [{\alpha}^{2}(z), (yx) \alpha (t)]$
\par
$ - (-1)^{\bar{z}(\bar{t} + \bar{x} + \bar{y}) + \bar{t}(\bar{x} + \bar{y})} [{\alpha}^{2}(z), as(x,t,y)]
+ (-1)^{\bar{z}(\bar{t} + \bar{x} + \bar{y}) + \bar{t}(\bar{x} + \bar{y})} [{\alpha}^{2}(z), (xt) \alpha (y)]$
\par
$ - (-1)^{\bar{z}(\bar{t} + \bar{x}) + \bar{y}(\bar{t} + \bar{z})} [{\alpha}^{2}(z), as(t,x,y)]
+ (-1)^{\bar{z}(\bar{t} + \bar{x}) + \bar{y}(\bar{t} + \bar{z})} [{\alpha}^{2}(z), (tx) \alpha (y)]$
\par
$ - (-1)^{\bar{z}(\bar{x} + \bar{y}) + \bar{t}(\bar{x} + \bar{z})} [{\alpha}^{2}(z), as(x,y,t)]
+ (-1)^{\bar{z}(\bar{x} + \bar{y}) + \bar{t}(\bar{x} + \bar{z})} [{\alpha}^{2}(z), (xy) \alpha (t)]$
\par
$ - (-1)^{\bar{z}(\bar{y} + \bar{t}) + \bar{x}(\bar{y} + \bar{z})} [{\alpha}^{2}(z), as(y,t,x)]
+ (-1)^{\bar{z}(\bar{y} + \bar{t}) + \bar{x}(\bar{y} + \bar{z})} [{\alpha}^{2}(z), (yt) \alpha (x)]$
\par
$ - (-1)^{\bar{z}(\bar{t} + \bar{x} + \bar{y}) + \bar{x}(\bar{t} + \bar{y})} [{\alpha}^{2}(z), as(y,x,t)]
- (-1)^{\bar{z}(\bar{t} + \bar{x} + \bar{y}) + \bar{x}(\bar{t} + \bar{y})} [{\alpha}^{2}(z), \alpha (y) (xt)]$
\par
$ - (-1)^{\bar{z}(\bar{t} + \bar{x} + \bar{y}) + \bar{t}(\bar{x} + \bar{y})} [{\alpha}^{2}(z), as(x,t,y)]
- (-1)^{\bar{z}(\bar{t} + \bar{x} + \bar{y}) + \bar{t}(\bar{x} + \bar{y})} [{\alpha}^{2}(z), \alpha (x) (ty)]$
\par
$ - (-1)^{\bar{z}(\bar{t} + \bar{x} + \bar{y}) + \bar{y}(\bar{t} + \bar{x})} [{\alpha}^{2}(z), as(t,y,x)]
- (-1)^{\bar{z}(\bar{t} + \bar{x} + \bar{y}) + \bar{y}(\bar{t} + \bar{x})} [{\alpha}^{2}(z), \alpha (t) (yx)]$
\par
$ - (-1)^{\bar{x}(\bar{t} + \bar{z}) + \bar{z}(\bar{t} + \bar{y})} [{\alpha}^{2}(z), as(x,y,t)]
- (-1)^{\bar{x}(\bar{t} + \bar{z}) + \bar{z}(\bar{t} + \bar{y})} [{\alpha}^{2}(z), \alpha (x) (yt)]$
\par
$ - (-1)^{\bar{y}(\bar{x} + \bar{z}) + \bar{z}(\bar{x} + \bar{t})} [{\alpha}^{2}(z), as(y,t,x)]
- (-1)^{\bar{y}(\bar{x} + \bar{z}) + \bar{z}(\bar{x} + \bar{t})} [{\alpha}^{2}(z), \alpha (y) (tx)]$
\par
$ - (-1)^{\bar{t}(\bar{y} + \bar{z}) + \bar{z}(\bar{x} + \bar{y})} [{\alpha}^{2}(z), as(t,x,y)]
- (-1)^{\bar{t}(\bar{y} + \bar{z}) + \bar{z}(\bar{x} + \bar{y})} [{\alpha}^{2}(z), \alpha (t) (xy)]$\\
and so, after rearranging terms, \\
\\
(4.5) \; ${\circlearrowleft}_{x,y,t} \; (-1)^{\bar{x}(\bar{y} + \bar{z}) + \bar{y}(\bar{t} + \bar{z})} 
\{ (-1)^{\bar{t}\bar{z}}{\alpha}^{2}(z)(\alpha (t)(yx))$ 
\par
$+(-1)^{\bar{t}\bar{z}+\bar{x}\bar{y}}{\alpha}^{2}(z)(\alpha (t)(xy))
+ (-1)^{(\bar{t}+\bar{z})(\bar{x}+\bar{y})} ((yx)\alpha (t)) {\alpha}^{2}(z)$
\par
$+ (-1)^{\bar{x}\bar{y} + (\bar{t}+\bar{z})(\bar{x}+\bar{y})} ((xy)\alpha (t)) {\alpha}^{2}(z)$
\par
$- (-1)^{\bar{z}(\bar{x} + \bar{y})} (\alpha (t) (yx)) {\alpha}^{2}(z)
- (-1)^{\bar{x}\bar{y} + \bar{z}(\bar{x} + \bar{y})} (\alpha (t) (xy)) {\alpha}^{2}(z)$
\par
$- (-1)^{\bar{t}(\bar{x} + \bar{y} + \bar{z})} {\alpha}^{2}(z) ((yx)\alpha (t))
- (-1)^{\bar{x}\bar{y} + \bar{t}(\bar{x} + \bar{y} + \bar{z})} {\alpha}^{2}(z) ((xy)\alpha (t)) \}$
\par
$= - (-1)^{\bar{z}(\bar{t} + \bar{x} + \bar{y}) + \bar{y}(\bar{t} + \bar{x})} [{\alpha}^{2}(z), as(t,y,x)]$
\par
$- (-1)^{\bar{z}(\bar{t} + \bar{x} + \bar{y}) + \bar{x}(\bar{t} + \bar{y})} [{\alpha}^{2}(z), as(y,x,t)]$
\par
$- (-1)^{\bar{z}(\bar{t} + \bar{x} + \bar{y}) + \bar{t}(\bar{x} + \bar{y})} [{\alpha}^{2}(z), as(x,t,y)]
- (-1)^{\bar{z}(\bar{t} + \bar{x}) + \bar{y}(\bar{t} + \bar{z})} [{\alpha}^{2}(z), as(t,x,y)]$
\par
$- (-1)^{\bar{z}(\bar{x} + \bar{y}) + \bar{t}(\bar{x} + \bar{z})} [{\alpha}^{2}(z), as(x,y,t)]
- (-1)^{\bar{z}(\bar{y} + \bar{t}) + \bar{x}(\bar{y} + \bar{z})} [{\alpha}^{2}(z), as(y,t,x)]$. \\
\\
Puting (4.5) in (4.4) we get (4.3), which completes the proof. \hfill $\square$  
\\
\par
Using (4.3) we now prove the main result of this section. \\

\par
{\bf Theorem 4.1.} {\it Let $\mathcal{A} = (A,*,\alpha)$ be a multiplicative right Hom-alternative superalgebra. Then
$\mathcal{A}$ is Hom-Jordan-admissible, i.e. ${\mathcal{A}}^{+}$ is a Hom-Jordan superalgebra}.\\

\par
{\it Proof.} We need to prove that (4.2) holds. Applying the right superalternativity (2.1) to the right-hand side
of (4.3), we get (for simplicity we replaced the product ''$*$`` by juxtaposition in suitable places)
\par
$8 {\circlearrowleft}_{x,y,t} \; (-1)^{\bar{x}(\bar{y} + \bar{z}) + \bar{y}(\bar{t} + \bar{z})} 
\{ (\alpha (t) \circ \alpha (z)) \circ \alpha (y \circ x)-{\alpha}^{2} (t) \circ ( \alpha (z) \circ (y \circ x)) \}$
\par
$= {\circlearrowleft}_{x,y,t} \; (-1)^{\bar{x}(\bar{y} + \bar{z}) + \bar{y}(\bar{t} + \bar{z})} 
\{ (-1)^{\bar{t}\bar{z}} as_{\mathcal{A}} (\alpha (z), \alpha (t), yx)$ 
\par
$+ (-1)^{\bar{t}\bar{z} + \bar{x}\bar{y}} as_{\mathcal{A}} (\alpha (z), \alpha (t), xy)$
\par
$ - (-1)^{\bar{t}(\bar{x} + \bar{y} + \bar{z})} as_{\mathcal{A}} (\alpha (z),yx,\alpha (t))
- (-1)^{\bar{x}\bar{y} + \bar{t}(\bar{x} + \bar{y} + \bar{z})} as_{\mathcal{A}} (\alpha (z),xy,\alpha (t)) \}$\\
i.e., again by (2.1) and the definition of the super-Jordan product ''$\circ$``,\\
\\
(4.6) \; $8 {\circlearrowleft}_{x,y,t} \; (-1)^{\bar{x}(\bar{y} + \bar{z}) + \bar{y}(\bar{t} + \bar{z})} 
\{ (\alpha (t) \circ \alpha (z)) \circ \alpha (y \circ x)-{\alpha}^{2} (t) \circ ( \alpha (z) \circ (y \circ x)) \}$
\par
\hspace{1.0cm}$= 4{\circlearrowleft}_{x,y,t} \; \{ (-1)^{\bar{t}(\bar{y} + \bar{z}) + \bar{z}(\bar{x} + \bar{y})}
as_{\mathcal{A}} (\alpha (z), \alpha (t), x \circ y) \}$. \\
\par
The following identity is proved to hold in any right Hom-alternative superalgebra (\cite{Iss2}, Theorem 3.6, identity
(3.14)): \\
\\
(4.7) \; $ as_{\mathcal{A}} (\alpha (z), \alpha (t), x \circ y) = (-1)^{\bar{x}\bar{y}} as_{\mathcal{A}}(z,t,y){\alpha}^{2}(x)
- (-1)^{\bar{t}\bar{x}} as_{\mathcal{A}}(z,x,t){\alpha}^{2}(y)$
\par
\hspace{4cm}$+ as_{\mathcal{A}} (\alpha (z),[t,x],\alpha (y)) + (-1)^{\bar{x}\bar{y}} 
as_{\mathcal{A}} (\alpha (z),[t,y],\alpha (x))$ \\
\\
(note that (4.7) above and (3.14) of \cite{Iss2} differ by coefficients because of specific definitions of ''$\circ$``
and ''$[\cdot,\cdot]$``). Therefore, replacing $ as_{\mathcal{A}} (\alpha (z), \alpha (t), x \circ y)$ in (4.6) with its
expression from (4.7) and next developing the sum ${\circlearrowleft}_{x,y,t}$, we obtain
\par
$8 {\circlearrowleft}_{x,y,t} \; (-1)^{\bar{x}(\bar{y} + \bar{z}) + \bar{y}(\bar{t} + \bar{z})} 
\{ (\alpha (t) \circ \alpha (z)) \circ \alpha (y \circ x)-{\alpha}^{2} (t) \circ ( \alpha (z) \circ (y \circ x)) \}$
\par
$= 4\{ (-1)^{\bar{x}\bar{y} + \bar{t}(\bar{y} + \bar{z}) + \bar{z}(\bar{x} + \bar{y})} as_{\mathcal{A}} (z,t,y){\alpha}^{2}(x)$
\par
$-(-1)^{\bar{t}\bar{x}+\bar{t}(\bar{y}+\bar{z}) + \bar{z}(\bar{x}+\bar{y})} as_{\mathcal{A}} (z,x,t){\alpha}^{2}(y)$
\par
$+ (-1)^{\bar{t}(\bar{y} + \bar{z}) + \bar{z}(\bar{x} + \bar{y})} as_{\mathcal{A}} (\alpha (z),[t,x], \alpha (y))$
\par
$+ (-1)^{\bar{x}\bar{y} + \bar{t}(\bar{y} + \bar{z}) + \bar{z}(\bar{x} + \bar{y})} as_{\mathcal{A}} (\alpha (z),[t,y], \alpha (x))$
\par
$+ (-1)^{\bar{t}\bar{x}+\bar{y}(\bar{x}+\bar{z}) + \bar{z}(\bar{t}+\bar{x})} as_{\mathcal{A}} (z,y,x){\alpha}^{2}(t)
- (-1)^{\bar{y}(\bar{t} + \bar{x} + \bar{z}) + \bar{z}(\bar{t}+\bar{x})} as_{\mathcal{A}} (z,t,y){\alpha}^{2}(x)$
\par
$+ (-1)^{\bar{x}(\bar{t} + \bar{z}) + \bar{z}(\bar{y} + \bar{t})} as_{\mathcal{A}} (\alpha (z),[x,y], \alpha (t))$
\par
$+ (-1)^{\bar{t}\bar{x} + \bar{y}(\bar{x} + \bar{z}) + \bar{z}(\bar{t} + \bar{x})} as_{\mathcal{A}} (\alpha (z),[y,x], \alpha (t))$
\par
$+ (-1)^{\bar{t}\bar{y}+ \bar{x}(\bar{t} + \bar{z}) + \bar{z}(\bar{y}+\bar{t})} as_{\mathcal{A}} (z,x,t){\alpha}^{2}(y)
- (-1)^{\bar{x}(\bar{t} + \bar{y} + \bar{z}) + \bar{z}(\bar{y}+\bar{t})} as_{\mathcal{A}} (z,y,x){\alpha}^{2}(t)$
\par
$+ (-1)^{\bar{y}(\bar{x} + \bar{z}) + \bar{z}(\bar{t} + \bar{x})} as_{\mathcal{A}} (\alpha (z),[y,t], \alpha (x))$
\par
$+ (-1)^{\bar{t}\bar{y} + \bar{x}(\bar{t} + \bar{z}) + \bar{z}(\bar{t} + \bar{y})} as_{\mathcal{A}} (\alpha (z),[x,t], \alpha (y)) \}$
\par
$= 0$ (by the right superalternativity and the superskewsymmetry of ''$[\cdot,\cdot]$``).\\
This completes the proof. \hfill $\square$  
\\
\par
As an immediate consequence, we have the following corollary (see also in \cite[Theorem 6.1]{AAM} but the proof uses
Lemma 6.1 therein).\\

\par
{\bf Corollary 4.1.} {\it Every multiplicative Hom-alternative superalgebra is Hom-Jordan- admissible}.\\

\par
{\it Proof.} Since a Hom-alternative superalgebra is first right Hom-alternative, a proof follows from Theorem 4.1. \hfill $\square$  

\section{Hom-Bol superalgebras. Construction theorems and example}

In \cite{AI1} Hom-Bol algebras were defined as a Hom-type generalization of Bol algebras. In this section we define
Hom-Bol superalgebras as a generalization both of Bol superalgebras \cite{Ruk}
and Hom-Bol algebras \cite{AI1}. Next we point out  some construction theorems.\\
\par
{\bf Definition 5.1.} A {\it Hom-Bol superalgebra} is a quadruple \\ ${\mathcal{A}}_{\alpha} := (A, [\cdot , \cdot], \{ \cdot, \cdot, \cdot \}, \alpha)$
where $A$ is a superspace, $[\cdot , \cdot]$ (resp. $\{ \cdot, \cdot, \cdot \}$) is a binary (resp. ternary) 
operation on $A$ such that $[A_{i},A_{j}] \subseteq A_{i+j}$, $\{ A_{i},A_{j},A_{k} \} \subseteq A_{i+j+k}$ and\\
\par
(SHB01) $\alpha ([x,y]) = [\alpha (x), \alpha (y)]$,
\par
(SHB02) $\alpha (\{x,y,z\}) = \{ \alpha (x), \alpha (y), \alpha (z) \}$,
\par
(SHB1) $[x,y] = - (-1)^{\bar{x} \bar{y}} [y,x]$,
\par
(SHB2) $\{x,y,z\} = - (-1)^{\bar{x} \bar{y}} \{y,x,z\}$,
\par
(SHB3) $\{x,y,z\} + (-1)^{\bar{x} (\bar{y} + \bar{z})} \{y,z,x\} + (-1)^{\bar{z} (\bar{x} + \bar{y})}\{z,x,y\} = 0$,
\par
(SHB4) $\{ \alpha (x), \alpha (y), [u,v] \} = [\{x,y,u\},{{\alpha}^{2}}(v)] + (-1)^{\bar{u} (\bar{x} + \bar{y})}[{{\alpha}^{2}}(u),\{x,y,v\}]$
\par
\hspace{1.5cm}$+ (-1)^{(\bar{x} + \bar{y}) (\bar{u} + \bar{v})} (\{\alpha(u),\alpha(v),[x,y]\} - [[\alpha(u),\alpha(v)],[\alpha(x),\alpha(y)]])$,
\par
(SHB5) $\{{{\alpha}^{2}}(x),{{\alpha}^{2}}(y),\{u,v,w\} \} = \{ \{x,y,u\}, {{\alpha}^{2}}(v),{{\alpha}^{2}}(w)\}$
\par
\hspace{5.9cm}$+ (-1)^{\bar{u} (\bar{x} + \bar{y})} \{{{\alpha}^{2}}(u), \{x,y,v\}, {{\alpha}^{2}}(w)\}$
\par
\hspace{5.9cm}$+ (-1)^{(\bar{x} + \bar{y}) (\bar{u} + \bar{v})} \{{{\alpha}^{2}}(u),{{\alpha}^{2}}(v), \{x,y,w\} \}$ \\
\\
for all homogeneous $u,v,w,x,y,z \in A$.\\
\par
We observe that for $\alpha = Id$, any Hom-Bol superalgebra reduces to a Bol superalgebra and a Hom-Bol superalgebra with a zero odd part is a 
Hom-Bol algebra. If $[x,y]=0$ for all homogeneous $x,y \in A$, one gets a {\it Hom-Lie supertriple system} 
$(A,\{ \cdot, \cdot, \cdot \}, {\alpha}^2)$. We note that, for Hom-Lie supertriple systems, construction theorems could
be proved in full analogy as for Hom-Lie triple systems (see \cite[Theorem 3.3, Corollary 3.4 and Corollary 3.5]{Yau5}).\\
\par

{\bf Theorem 5.1.} {\it Let ${\mathcal{A}}_{\alpha} := (A, [\cdot , \cdot], \{ \cdot, \cdot, \cdot \}, \alpha)$ be a Hom-Bol superalgebra and 
$\beta : A \rightarrow A$ an even self-morphism of ${\mathcal{A}}_{\alpha}$ i.e. $\beta \circ [\cdot, \cdot]
= [\cdot, \cdot] \circ {\beta}^{\otimes 2}$, $\beta \circ \{ \cdot, \cdot, \cdot \}
= \{ \cdot, \cdot, \cdot \} \circ {\beta}^{\otimes 3}$, and $\alpha \circ \beta = \beta \circ \alpha$. Let ${\beta}^0 = Id$
and ${\beta}^n = \beta \circ {\beta}^{n-1}$ for any integer $n \geq 0$ and define on $A$ a binary operation $[\cdot, \cdot]_{{\beta}^n}$ and a
ternary operation $\{ \cdot, \cdot, \cdot \}_{{\beta}^n}$ by
\par
$[x, y]_{{\beta}^n} := {\beta}^n ([x,y])$,
\par
$\{ x,y,z \}_{{\beta}^n} := {\beta}^{2n} (\{x,y,z\})$. \\
Then ${\mathcal{A}}_{{\beta}^n} := (A, [\cdot, \cdot]_{{\beta}^n}, \{ \cdot, \cdot, \cdot \}_{{\beta}^n}, {\beta}^n \circ \alpha)$ is a Hom-Bol
superalgebra}. \\

\par
{\it Proof.} The proof is similar to that of \cite[Theorem 3.2]{AI1}. \hfill $\square$ \\ 
\par
From Theorem 5.1 we get the following extension of the Yau twisting principle \cite{Yau1} giving a construction of Hom-Bol superalgebras from
Bol superalgebras.\\
\par

{\bf Corollary 5.1.} {\it Let $(A, [\cdot , \cdot], \{ \cdot, \cdot, \cdot \})$ be a Bol superalgebra and $\beta$ an even self-morphism of
$(A, [\cdot , \cdot], \{ \cdot, \cdot, \cdot \})$. Define on $A$ a binary operation $[\cdot, \cdot]_{\beta}$ and a ternary operation 
$\{ \cdot, \cdot, \cdot \}_{\beta}$ by
\par
$[x, y]_{\beta} := {\beta} ([x,y])$,
\par
$\{ x,y,z \}_{\beta} := {\beta}^{2} (\{x,y,z\})$. \\
Then ${\mathcal{A}}_{\beta} := (A, [\cdot, \cdot]_{\beta}, \{ \cdot, \cdot, \cdot \}_{\beta}, \beta )$ is a Hom-Bol superalgebra. 
Moreover, if \\ $(A', [\cdot , \cdot]', \{ \cdot, \cdot, \cdot \}')$ is another Bol superalgebra, ${\beta}'$ 
an even self-morphism of \\ $(A', [\cdot , \cdot]', \{ \cdot, \cdot, \cdot \}')$
and if $f:A \rightarrow A'$ is a Bol superalgebra even morphism satisfying $f \circ \beta = {\beta}' \circ f$, then 
$f:{\mathcal A}_{\beta} \rightarrow {\mathcal{A}'}_{{\beta}'}$ is a morphism of Hom-Bol superalgebras, where 
${\mathcal{A}'}_{{\beta}'} := (A', [\cdot , \cdot{]'}_{{\beta}'}, \{ \cdot, \cdot, \cdot {\}'}_{{\beta}'}, {\beta}')$}.\\
\par
{\it Proof.} The first part of the corollary comes from Theorem 5.1 when $n=1$. The second part is proved in a similar way as 
\cite[Corollary 4.5]{Iss1}. \hfill $\square$ \\ 

\par 
{\bf Example 5.1.} Let  $A = A_0 \oplus A_1$ be a superspace over a field of characteristic not $2$ where $A_0$ is a $1$-dimensional vector space 
with basis $\{i\}$ and  $A_1$ a $2$-dimensional vector space with basis $\{j,k\}$. Define on $A$ the following only nonzero products on basis elements:
\par
$i*j=k$;
\par
$j*i=k$, $j*k=2i$;
\par
$k*j=4i$.\\
Then $(A, *)$ is a right alternative superalgebra \cite{Sh2}. Now consider on $(A, *)$ the supercommutator $[\cdot , \cdot]$ and the ternary operation
defined as
\par
$\{x,y,z\} := (-1)^{\bar{x}(\bar{y} + \bar{z})} as_{\mathcal{A}^{+}}(y,z,x)$.\\
Then Theorem 3.1 implies that $(A,[\cdot,\cdot],\{ \cdot,\cdot,\cdot \})$ is a Bol superalgebra, where the only nonzero products are:
\par
$[j,k] = 6i$;
\par
$[k,j] = 6i$;
\par
$\{i,j,j\} = 4i$;
\par
$\{j,i,j\} = -4i$, $\{j,j,i\} = -8i$, $\{j,j,k\} = -8k$, $\{j,k,j\} = 4k$;
\par
$\{k,j,j\} = 4k$.\\
Next define a linear map $\beta : A \rightarrow A$ by setting
\par
$\beta (i) = ai$, $\beta (j) = j+bk$, $\beta (k) = ak$\\
with $a \neq 0$. Then it is easily seen that $\beta$ is an even self-morphism of \\ $(A,[\cdot,\cdot],\{ \cdot,\cdot,\cdot \})$ and Corollary 5.1 implies
that ${\mathcal{A}}_{\beta} := (A, [\cdot, \cdot]_{\beta}, \{ \cdot, \cdot, \cdot \}_{\beta}, \beta )$ is a Hom-Bol superalgebra with products given as
\par
$[j,k]_{\beta} = 6ai$;
\par
$[k,j]_{\beta} = 6ai$;
\par
$\{i,j,j\}_{\beta} = 4a^{2}i$;
\par
$\{j,i,j\}_{\beta} = -4a^{2}i$, $\{j,j,i\}_{\beta} = -8a^{2}i$, $\{j,j,k\}_{\beta} = -8a^{2}k$, $\{j,k,j\}_{\beta} = 4a^{2}k$;
\par
$\{k,j,j\}_{\beta} = 4a^{2}k$.\\
\par
The notion of an $n$th derived (binary) Hom-algebra of a given Hom-algebra is first introduced in \cite{Yau4} and the closure of a given type of 
Hom-algebras under taking $n$th derived
Hom-algebras is a property that is characteristic of the variety of Hom-algebras. Later on, this notion is extended to binary-ternary Hom-algebras \cite{AI1} 
(for binary-ternary Hom-superalgebras \cite{GI}, the notion is the same as in the case of binary-ternary Hom-algebras).\\
\par
{\bf Definition 5.2.}  \cite{GI} Let ${\mathcal A} := (A,*,\{ \cdot,\cdot,\cdot \},\alpha)$ be a binary-ternary Hom-superalgebra and $n\geq 0$ an integer. 
Define on $A$ the $n$th derived 
binary operation $*^{(n)}$ and the $n$th derived ternary operation $\{ \cdot,\cdot,\cdot \}^{(n)}$ by
\par
$x *^{(n)} y := \alpha^{{2^n}-1}(x * y)$,
\par
$\{x,y,z\}^{(n)} := \alpha^{{2^{n+1}}-2}(\{x,y,z\})$, \\
for all homogeneous $x,y,z$ in A. Then ${\mathcal A}^{(n)} := (A, *^{(n)},\{ \cdot,\cdot,\cdot \}^{(n)}, \alpha^{2^{n}})$ is called the $n$th 
{\it derived  (binary-ternary) 
Hom-superalgebra} of $\mathcal A$. \\
\par 
As for Hom-Bol algebras, the category of Hom-Bol superalgebras is closed under taking derived Hom-superalgebras as stated in the following\\
\par

{\bf Theorem 5.2.} {\it Let ${\mathcal A} := (A,[\cdot,\cdot],\{ \cdot,\cdot,\cdot \},\alpha)$ be a Hom-Bol superalgebra.
Then, for each $n \geq 0$, the $n$th derived Hom-superalgebra ${\mathcal A}^{(n)}$ is a Hom-Bol superalgebra. In
particular, the $n$th derived Hom-superalgebra of a Hom-Lie supertriple system is a Hom-Lie supertriple system}. \\

\par
{\it Proof.} The proof of the first part of the theorem is similar to that of \cite[Theorem 3.5]{AI1}, and the second 
part follows immediately. \hfill $\square$  

\section{Jordan and Hom-Jordan supertriple systems}

The theory of Jordan triple systems is developed in \cite{Jac, Loos, Meyb} and a Hom-type generalization
of Jordan triple systems (called Hom-Jordan triple systems) is considered in \cite{Yau5}. As far as we could find from
existing literature, the $\mathbb{Z}_2$-graded generalization of Jordan triple systems (call them Jordan supertriple
systems) is first considered in \cite{BC}. In this section we give an example of a Jordan supertriple system and next 
consider a Hom-type generalization of Jordan supertriple systems.\\
\par
{\bf Definition 6.1.} \cite{BC} A {\it Jordan supertriple system} is a pair $(V, \langle \cdot , \cdot , \cdot \rangle)$
consisting of a superspace $V = V_0 \oplus V_1$ and a trilinear operation $\langle \cdot , \cdot , \cdot \rangle : V \times
V \times V \rightarrow V$ such that 
\par
$\bullet$ $\langle V_i , V_j , V_k \rangle \subseteq V_{i+j+k}$,
\par
$\bullet$ $\langle x,y,z \rangle = (-1)^{\bar{x}\bar{y}+\bar{x}\bar{z}+\bar{y}\bar{z}} \langle z,y,x \rangle$ ({\it outer
supersymmetry}),
\par
$\bullet$ $\langle x,y, \langle u,v,w \rangle \rangle - \langle \langle x,y,u \rangle, v,w \rangle =
(-1)^{(\bar{x}+\bar{y})(\bar{u}+\bar{v})} (\langle u,v, \langle x,y,w \rangle \rangle$ 
\par
\hspace{9.0cm}$- \langle u, \langle v,x,y \rangle , w \rangle)$ \\ ({\it Jordan supertriple identity})\\
for all homogeneous $u,v,w,x,y,z$ in $V$. \\
\par
The following example is inspired by a construction of Jordan triple systems associated to symmetric bilinear forms
on any finite-dimensional vector space over a field of characteristic zero \cite{Loos}. \\
\par
{\bf Example 6.1.} Let $V$ be a superspace over a field $\mathbb{K}$ of characteristic zero and $\langle \cdot | \cdot 
\rangle : V \times V \rightarrow \mathbb{K}$ be a supersymmetric bilinear form on $V$ i.e. $\langle x | y \rangle = 
(-1)^{\bar{x}\bar{y}} \langle y | x \rangle$ for all homogeneous $x,y \in V$. Consider on $V$ the trilinear operation
$\langle \cdot , \cdot , \cdot \rangle$ defined by
\par
$\langle x,y,z \rangle := \lambda ( \langle x | y \rangle z + (-1)^{\bar{x}(\bar{y}+\bar{z})} \langle y | z \rangle x
-(-1)^{\bar{z}(\bar{x}+\bar{y})} \langle z | x \rangle y)$,\\
where $\lambda \in \mathbb{K}$. We have
\par
$(-1)^{\bar{x}\bar{y}+\bar{x}\bar{z}+\bar{y}\bar{z}} \langle z,y,x \rangle = 
\lambda (-1)^{\bar{x}\bar{y}+\bar{x}\bar{z}+\bar{y}\bar{z}} ( \langle z|y \rangle x + (-1)^{\bar{z}(\bar{y}+\bar{x})} 
\langle y | x \rangle z$ 
\par
$-(-1)^{\bar{x}(\bar{z}+\bar{y})} \langle x | z \rangle y)$
\par
$=\lambda (-1)^{\bar{x}\bar{y}+\bar{x}\bar{z}+\bar{y}\bar{z}} ( (-1)^{\bar{y}\bar{z}}\langle y|z \rangle x 
+ (-1)^{\bar{z}(\bar{y}+\bar{x})+ \bar{x}\bar{y}} 
\langle x | y \rangle z -(-1)^{\bar{x}\bar{y}} \langle z | x \rangle y)$
\par
$=\langle x,y,z \rangle$\\
i.e. we get the outer symmetry. Next, from one hand we have
\par
$\langle x,y, \langle u,v,w \rangle \rangle - \langle \langle x,y,u \rangle, v,w \rangle$
\par
$= \langle x,y, \lambda ( \langle u | v \rangle w + (-1)^{\bar{u}(\bar{v}+\bar{w})} \langle v | w \rangle u -(-1)^{\bar{w}(\bar{u}+\bar{v})} \langle w | u \rangle v)$
\par
$- \langle \lambda ( \langle x | y \rangle u + (-1)^{\bar{x}(\bar{y}+\bar{u})} \langle y | u \rangle x
-(-1)^{\bar{u}(\bar{x}+\bar{y})} \langle u | x \rangle y), v, w \rangle$
\par
$= (-1)^{(\bar{x}+\bar{y})(\bar{u}+\bar{v})+\bar{x}( \bar{y} + \bar{w})} {\lambda}^{2} \langle u | v \rangle
\langle y | w \rangle x - (-1)^{(\bar{x}+\bar{y})(\bar{u}+\bar{v})+\bar{w}( \bar{x} + \bar{y})} {\lambda}^{2} \langle u | v \rangle
\langle w | x \rangle y$
\par
$- (-1)^{\bar{w}(\bar{u}+\bar{v})+(\bar{x}+\bar{y})(\bar{u}+\bar{w})+\bar{x}(\bar{y}+\bar{v})} 
{\lambda}^{2} \langle w | u \rangle \langle y | v \rangle x$
\par
$+ (-1)^{\bar{w}(\bar{u}+\bar{v})+(\bar{x}+\bar{y})(\bar{u}+\bar{w})+\bar{v}(\bar{x}+\bar{y})} 
{\lambda}^{2} \langle w | u \rangle \langle v | x \rangle y$
\par
$- (-1)^{\bar{x}(\bar{y}+\bar{u})} {\lambda}^{2} \langle y | u \rangle \langle x | v \rangle w
+ (-1)^{\bar{x}(\bar{y}+\bar{u}) + \bar{w}(\bar{x}+\bar{v})} {\lambda}^{2} \langle y | u \rangle \langle w | x \rangle v$
\par
$+(-1)^{\bar{u}(\bar{x}+\bar{y})} {\lambda}^{2} \langle u | x \rangle \langle y | v \rangle w 
- (-1)^{\bar{u}(\bar{x}+\bar{y}) + \bar{w}(\bar{y}+\bar{v})} {\lambda}^{2} \langle u | x \rangle \langle w | y \rangle v$\\
and, from the other hand,
\par
$(-1)^{(\bar{x}+\bar{y})(\bar{u}+\bar{v})} (\langle u,v, \langle x,y,w \rangle \rangle - \langle u, \langle v,x,y 
\rangle , w \rangle) $
\par
$= (-1)^{(\bar{x}+\bar{y})(\bar{u}+\bar{v})} \langle u,v, \lambda ( \langle x | y \rangle w + (-1)^{\bar{x}(\bar{y}+\bar{w})} 
\langle y | w \rangle x -(-1)^{\bar{w}(\bar{x}+\bar{y})} \langle w | x \rangle y) \rangle$
\par
$-(-1)^{(\bar{x}+\bar{y})(\bar{u}+\bar{v})} \langle u, \lambda ( \langle v | x \rangle y 
+ (-1)^{\bar{v}(\bar{x}+\bar{y})} 
\langle x | y \rangle v -(-1)^{\bar{y}(\bar{v}+\bar{x})} \langle y | v \rangle x), w \rangle$
\par
$= (-1)^{(\bar{x}+\bar{y})(\bar{u}+\bar{v})+\bar{x}( \bar{y} + \bar{w})} {\lambda}^{2} \langle u | v \rangle
\langle y | w \rangle x - (-1)^{\bar{u}\bar{y} + \bar{v}\bar{w}} 
{\lambda}^{2} \langle x | u \rangle \langle y | w \rangle v$
\par
$- (-1)^{(\bar{x}+\bar{y})(\bar{u}+\bar{v}+\bar{w})} {\lambda}^{2} \langle u | v \rangle \langle w | x \rangle y
+ (-1)^{\bar{x}(\bar{y}+\bar{u}) + \bar{w}(\bar{x}+\bar{v})} {\lambda}^{2} \langle y|u \rangle \langle w|x \rangle v$
\par
$- (-1)^{\bar{x}(\bar{u}+\bar{v}) + \bar{y}(\bar{u}+\bar{x})} {\lambda}^{2} \langle u|y \rangle \langle v|x \rangle w
+ (-1)^{\bar{v}(\bar{x}+\bar{y})+ \bar{w}(\bar{u}+\bar{v})+(\bar{u}+\bar{w})(\bar{x}+\bar{y})} 
{\lambda}^{2} \langle w | u \rangle \langle v | x \rangle y$
\par
$+ (-1)^{\bar{u}(\bar{x}+\bar{y})} {\lambda}^{2} \langle u | x \rangle \langle y | v \rangle w 
- (-1)^{\bar{x}(\bar{u}+\bar{v}+\bar{w}+\bar{y})+\bar{u}(\bar{w}+\bar{y}) + \bar{w}(\bar{v}+\bar{y})} 
{\lambda}^{2} \langle w | u \rangle \langle y | v \rangle x$. \\
Therefore the Jordan supertriple identity is verified and so $(V, \langle \cdot , \cdot , \cdot \rangle)$ is a Jordan
supertriple system.\\
\par
A Hom-version of a Jordan supertriple system is defined as follows.\\
\par
{\bf Definition 6.2.} A (multiplicative) {\it Hom-Jordan supertriple system} is a triple \\
$(V, \langle \cdot , \cdot , \cdot \rangle, \theta)$
consisting of a superspace $V = V_0 \oplus V_1$, a trilinear operation \\
$\langle \cdot , \cdot , \cdot \rangle : V \times V \times V \rightarrow V$, and a linear even self-map of $V$ such that 
\par
$\bullet$ $\langle V_i , V_j , V_k \rangle \subseteq V_{i+j+k}$,
\par
$\bullet$ $\langle x,y,z \rangle = (-1)^{\bar{x}\bar{y}+\bar{x}\bar{z}+\bar{y}\bar{z}} \langle z,y,x \rangle$ ({\it outer
supersymmetry}),
\par
$\bullet$ $\langle \theta (x), \theta (y), \langle u,v,w \rangle \rangle 
- \langle \langle x,y,u \rangle, \theta (v), \theta (w) \rangle $
\par
\hspace{0.3cm} $= (-1)^{(\bar{x}+\bar{y})(\bar{u}+\bar{v})} (\langle \theta (u), \theta (v), \langle x,y,w \rangle \rangle 
- \langle \theta (u), \langle v,x,y \rangle , \theta (w) \rangle)$ ({\it Hom-Jordan supertriple identity})\\
for all homogeneous $u,v,w,x,y,z$ in $V$. \\
\par
One observes that $V_0$ is nothing but a (multiplicative) {\it Hom-Jordan triple system} that is defined in \cite{Yau5}
while, for $\theta = Id$, $(V, \langle \cdot , \cdot , \cdot \rangle)$ is a Jordan supertriple system. We have the 
following construction theorem. \\
\par

{\bf Theorem 6.1.} {\it Let $V_{\alpha} := (V, \langle \cdot , \cdot , \cdot \rangle, \alpha)$ be a multiplicative 
Hom-Jordan supertriple system and $\beta : V \rightarrow V$ be an even self-morphism of $V_{\alpha}$ (i.e. 
$\beta \circ \langle \cdot , \cdot , \cdot \rangle = \langle \cdot , \cdot , \cdot \rangle \circ {\beta}^{\otimes 3}$
and $\beta \circ \alpha = \alpha \circ \beta$). Let ${\beta}^{0} = Id$ and ${\beta}^{n} = \beta \circ {\beta}^{n-1}$
for any integer $n \geq 0 $. Then the ternary Hom-superalgebra $V_{{\beta}^{n}} := (V, \langle \cdot , \cdot , \cdot 
{\rangle}_{{\beta}^{n}} := {\beta}^{n} \circ \langle \cdot , \cdot , \cdot \rangle, {\beta}^{n} \circ \alpha)$ is a 
multiplicative Hom-Jordan supertriple system}. \\

\par
{\it Proof.} First observe that the multiplicativity of $V_{{\beta}^{n}}$ follows from the one of $V_{\alpha}$ and from 
that $\beta \circ \alpha = \alpha \circ \beta$. Next we have
\par
$\langle x , y , z {\rangle}_{{\beta}^{n}} = {\beta}^{n} (\langle x , y , z \rangle) 
= (-1)^{\bar{x}\bar{y}+\bar{x}\bar{z}+\bar{y}\bar{z}} {\beta}^{n} (\langle z , y , x \rangle)$
\par
\hspace{1.9cm}$= (-1)^{\bar{x}\bar{y}+\bar{x}\bar{z}+\bar{y}\bar{z}} \langle z , y , x {\rangle}_{{\beta}^{n}}$ \\
and 
\par
$\langle ( {\beta}^{n} \circ \alpha )(x) , ( {\beta}^{n} \circ \alpha )(y), \langle u,v,w {\rangle}_{{\beta}^{n}} {\rangle}_{{\beta}^{n}}
- \langle \langle x,y,u {\rangle}_{{\beta}^{n}} , ({\beta}^{n} \circ \alpha)(v), ({\beta}^{n} \circ \alpha)(w) {\rangle}_{{\beta}^{n}}$
\par
$= {\beta}^{2n} (\langle \alpha (x), \alpha (y), \langle u,v,w \rangle \rangle 
- \langle \langle x,y,u \rangle, \alpha (v), \alpha (w) \rangle)$
\par
$= {\beta}^{2n} ( (-1)^{(\bar{x}+\bar{y})(\bar{u}+\bar{v})} ( \langle \alpha (u), \alpha (v), \langle x,y,w \rangle \rangle 
- \langle \alpha (u), \langle v,x,y \rangle , \alpha (w) \rangle) )$
\par
$= (-1)^{(\bar{x}+\bar{y})(\bar{u}+\bar{v})} ( {\beta}^{n}(\langle ({\beta}^{n} \circ \alpha)(u), 
({\beta}^{n} \circ \alpha)(v), {\beta}^{n} (\langle x , y , w \rangle) \rangle)$
\par
$- {\beta}^{n} ( \langle ({\beta}^{n} \circ \alpha)(u), {\beta}^{n} (\langle v , x , y \rangle), ({\beta}^{n} \circ 
\alpha)(w) \rangle) )$
\par
$= (-1)^{(\bar{x}+\bar{y})(\bar{u}+\bar{v})} ( \langle ({\beta}^{n} \circ \alpha)(u), ({\beta}^{n} \circ \alpha)(v),
\langle x , y , w {\rangle}_{{\beta}^{n}} {\rangle}_{{\beta}^{n}}$
\par
$ - \langle ({\beta}^{n} \circ \alpha)(u),
\langle v , x , y {\rangle}_{{\beta}^{n}}, ({\beta}^{n} \circ \alpha)(w) {\rangle}_{{\beta}^{n}} )$ \\
for all $u,v,w,x,y,z \in V$. So we get that $V_{{\beta}^{n}}$ is a multiplicative Hom-Jordan supertriple system.
\hfill $\square$ \\ 
\par
Now the Yau twisting principle gives rise to examples of Hom-Jordan supertriple systems. \\
\par

{\bf Corollary 6.1.} {\it Let $(V, \langle \cdot , \cdot , \cdot \rangle)$ be a Jordan supertriple system and $\beta$
an even self-morphism of $(V, \langle \cdot , \cdot , \cdot \rangle)$. If define on $V$ a ternary operation
$\langle \cdot , \cdot , \cdot {\rangle}_{\beta}$ by 
$\langle x, y, z {\rangle}_{\beta} = \beta (\langle x, y , z \rangle )$ for all homogeneous
$x,y,z \in V$, then $(V, \langle \cdot , \cdot , \cdot {\rangle}_{\beta}, \beta)$ is a Hom-Jordan supertriple system}.\\

\par
{\it Proof.} A proof follows from Theorem 6.1 when $\alpha = Id$ and $n=1$. \hfill $\square$ \\ 
\par
{\bf Example 6.2.} Let $(V, \langle \cdot , \cdot , \cdot \rangle)$ be the Jordan supertriple system of Example 6.1 and
let $\beta : V \rightarrow V$ be any even linear map such that $\langle \beta (x) | \beta (y) \rangle = 
\langle x | y \rangle$ for all homogeneous $x,y \in V$. Then $\beta$ is a self-morphism of $(V, \langle \cdot , \cdot , \cdot \rangle)$
and Corollary 6.1 implies that $(V, \langle \cdot , \cdot , \cdot {\rangle}_{\beta} = \beta \circ \langle \cdot , \cdot , \cdot \rangle,
\beta)$ is a multiplicative Hom-Jordan supertriple system.

\section{Hom-Bol structures on right Hom-alternative superalgebras}

In this section we first prove that Hom-Jordan superalgebras give rise to Hom-Jordan supertriple systems (Theorem 7.1).
Using this fact and exploiting the connection between Hom-Jordan supertriple systems and Hom-Lie supertriple systems
(Theorem 7.2), we next prove our main result connecting right Hom-alternative superalgebras with Hom-Bol superalgebras
(Theorem 7.3).

\subsection{From Hom-Jordan superalgebras to Hom-Jordan supertriple systems} 
\par
Jordan algebras are known to give rise to Jordan triple systems. In \cite{Yau3} a Hom-version of this result is proved,
and its $\mathbb{Z}_2$-graded generalization is mentioned in \cite{BC}. In this subsection we prove a combined version
of results from \cite{Yau3} and \cite{BC}. In our proof we shall consider a $\mathbb{Z}_2$-graded generalization of the
method developed in the section 3 of \cite{Yau3}. Thus we shall prove the following\\
\par

{\bf Theorem 7.1.} {\it Let $(A, \circ , \alpha)$ be a multiplicative Hom-Jordan superalgebra. Define the ternary
operation} \\
\\
(7.1) \; $\langle x, y , z \rangle := (xy) \alpha (z) + \alpha (x) (yz) - (-1)^{\bar{x}\bar{y}} \alpha (y) (xz)$ \\
\\
{\it for all homogeneous $x,y,z \in A$, where $xy := x \circ y$. Then $(A, \langle \cdot , \cdot , \cdot \rangle, 
{\alpha}^{2})$ is a multiplicative Hom-Jordan supertriple system}.\\

\par
The proof of Theorem 7.1 requires some notions and lemmas which are a $\mathbb{Z}_2$-graded generalization of 
Definition 3.2 and Lemmas 3.3-3.10 in \cite{Yau3}. Applying the Kaplansky rule to Definition 3.2 in \cite{Yau3}, we
have the following \\
\par
{\bf Definition 7.1.} Let $(A, \cdot , \alpha)$ be any Hom-superalgebra. For homogeneous $w,x,y,z \in A$, define the 
linear maps $L(x)$, $L(x,y)$, $L(x,y,z)$, $L(w,x,y,z)$ and $L_{x,y}$ by
\par
$L(x)(w) = x \cdot w$,
\par
$L(x,y) = L(\alpha (x) ) L(y) - (-1)^{\bar{x}\bar{y}} L(\alpha (y) ) L(x)$,
\par
$L(x,y,z) = L(\alpha (x), \alpha (y) )L(z) - (-1)^{\bar{z}(\bar{x}+\bar{y})} L( {\alpha}^{2} (z) ) L(x,y)$,
\par
$L(w,x,y,z) = L( \alpha (w), \alpha (x), \alpha (y) )L(z) - (-1)^{\bar{z}(\bar{w}+\bar{x}+\bar{y})} 
L( {\alpha}^{3} (z) ) L(w,x,y)$,
\par
$L_{x,y} = L(x \cdot y) \alpha + L(x,y)$. \\
\par
The composition of maps is denoted by juxtaposition. Observe that the operator $L_{x,y}$ is closely related to the product (7.1). In fact we have
$L_{x,y} (z) = L(xy)( \alpha (z) ) + L(x,y)(z)$ and so \\
\\
(7.2) \; $L_{x,y} = \langle x, y , - \rangle$. \\
\\
\par
From Definition 7.1 we have the following obvious properties that we shall use in the rest of this subsection.\\
\par

{\bf Lemma 7.1.} {\it Let $(A, \cdot , \alpha)$ be a multiplicative supercommutative Hom- superalgebra. Then, for all
homogeneous $x,y,z \in A$, we have
\par
(i) $L(x)(y) = (-1)^{\bar{x}\bar{y}} L(y)(x)$;
\par
(ii) $L(x+y) = L(x) + L(y)$;
\par
(iii) $L(x,y) = (-1)^{\bar{x}\bar{y}} L(y,x)$;
\par
(iv) $\alpha L(x) = L(\alpha (x)) \alpha$;
\par
(v) $\alpha L(x,y) = L(\alpha (x), \alpha (y)) \alpha$;
\par
(vi) $\alpha L(x,y,z) = L(\alpha (x), \alpha (y), \alpha (z)) \alpha$}. \hfill $\square$ \\

\par 
One notes that for the properties (iv)-(vi), the multiplicativity of $(A, \cdot , \alpha)$ is essential. The Lemmas
7.2-7.8 below will be used in the proof of Theorem 7.1.\\
\par

{\bf Lemma 7.2.} {\it Let $(A, \circ , \alpha)$ be a multiplicative Hom-Jordan superalgebra. Then\\
\par
${\circlearrowleft}_{w,x,z}  (-1)^{\bar{w}(\bar{x} + \bar{z})} L(\alpha (x), wz) \alpha = 0$ \\
\\
for all homogeneous $w,x,z \in A$, where $uv := u \circ v$}. \\

\par
{\it Proof.} Write the Hom-Jordan superidentity (4.2) as \\ ${\circlearrowleft}_{w,x,z}  
(-1)^{\bar{z}(\bar{w} + \bar{y}) + \bar{w}(\bar{x} + \bar{y})} {as}_{{\mathcal{A}}^{+}} (\alpha (x), \alpha (y), wz) =0$.
Then
\par
$0 = {\circlearrowleft}_{w,x,z}  
(-1)^{\bar{z}(\bar{w} + \bar{y}) + \bar{w}(\bar{x} + \bar{y})} {as}_{{\mathcal{A}}^{+}} (\alpha (x), \alpha (y), wz)$
\par
\hspace{0.3cm}$ = {\circlearrowleft}_{w,x,z}  (-1)^{\bar{z}(\bar{w} + \bar{x})} \{ \alpha (wz) (\alpha (x) \alpha (y))
- (-1)^{\bar{x}(\bar{w} + \bar{z})} {\alpha}^{2} (x) ((wz) \alpha (y)) \}$
\par
\hspace{0.3cm}$ = {\circlearrowleft}_{w,x,z}  (-1)^{\bar{z}(\bar{w} + \bar{x})} \{ L(\alpha (wz)) L(\alpha (x))
- (-1)^{\bar{x}(\bar{w} + \bar{z})} L( {\alpha}^{2} (x)) L(wz)\} (\alpha (y))$
\par
\hspace{0.3cm}$ = {\circlearrowleft}_{w,x,z}  (-1)^{\bar{z}(\bar{w} + \bar{x})} L(wz,\alpha (x)) (\alpha (y))$
\par
\hspace{0.3cm}$ = {\circlearrowleft}_{w,x,z} - (-1)^{\bar{w}(\bar{x} + \bar{z})} L( \alpha (x), wz) (\alpha (y))$\\
which proves the lemma. \hfill $\square$ \\
\par

{\bf Lemma 7.3.} {\it Let $(A, \circ , \alpha)$ be a multiplicative Hom-Jordan superalgebra. Then}\\
\\
(7.3) \; $L(L(x,y)(z)) {\alpha}^{2} = L(x,y,z)$ \hfill \\
\\
{\it for all homogeneous $x,y,z \in A$}. \\

\par
{\it Proof.} First observe that, by definition, \\
\\
(7.4) \; $L(x,y,z) = L({\alpha}^{2} (x)) L(\alpha (y)) L(z) - (-1)^{\bar{x}\bar{y}} L( {\alpha}^{2} (y)) L(\alpha (x)) L(z)$
\par
$- (-1)^{\bar{z}(\bar{x} + \bar{y})} L({\alpha}^{2} (z)) L(\alpha (x)) L(y)
+ (-1)^{\bar{x}\bar{y} + \bar{z}(\bar{x} + \bar{y})} L({\alpha}^{2} (z)) L(\alpha (y)) L(x)$. \\
\\
Now write the Hom-Jordan superidentity (4.2) as
\par
$0=(-1)^{\bar{y}(\bar{w} + \bar{z}) + \bar{w}(\bar{x} + \bar{z})} \{ \alpha (xy) (\alpha (w)\alpha (z)) 
- {\alpha}^{2} (x) (\alpha (y) (wz)) \}$
\par
\hspace{0.3cm}$+ (-1)^{\bar{y}(\bar{z} + \bar{x}) + \bar{z}(\bar{w} + \bar{x})} \{ (\alpha (w)\alpha (y)) 
\alpha (zx) - {\alpha}^{2} (w) (\alpha (y) (zx)) \}$
\par
\hspace{0.3cm}$+ (-1)^{\bar{y}(\bar{x} + \bar{w}) + \bar{x}(\bar{z} + \bar{w})} \{ (\alpha (z)\alpha (y)) 
\alpha (xw) - {\alpha}^{2} (z) (\alpha (y) (xw)) \}$
\par
\hspace{0.3cm}$= (-1)^{\bar{y}(\bar{w} + \bar{z}) + \bar{w}(\bar{x} + \bar{z})} \{ (-1)^{\bar{w}\bar{z}} 
L( \alpha (xy)) L(\alpha (z)) \alpha$ 
\par
\hspace{4.0cm}$- (-1)^{\bar{w}\bar{z}} L({\alpha}^{2} (x)) L(\alpha (y)) L(z) \} (w)$
\par
\hspace{0.3cm}$+ (-1)^{\bar{y}(\bar{z} + \bar{x}) + \bar{z}(\bar{w} + \bar{x})} \{ (-1)^{\bar{w}\bar{z}+ 
(\bar{z} + \bar{x})(\bar{w} + \bar{y})} L( \alpha (zx)) L(\alpha (y)) \alpha $ 
\par
\hspace{3.7cm}$- (-1)^{\bar{w}(\bar{x} + \bar{y} + \bar{z})} L(\alpha (y)(zx) {\alpha}^{2} \} (w)$
\par
\hspace{0.3cm}$+ (-1)^{\bar{y}(\bar{x} + \bar{w}) + \bar{x}(\bar{z} + \bar{w})} \{ L(\alpha (zy))L(\alpha (x))\alpha
- L({\alpha}^{2} (z)) L(\alpha (y)) L(x) \} (w)$\\
and so \\
\\
(7.5) \; $(-1)^{\bar{w}\bar{x} + \bar{y}(\bar{w} + \bar{z})} L( \alpha (xy)) L(\alpha (z)) \alpha
- (-1)^{\bar{w}\bar{x} + \bar{y}(\bar{w} + \bar{z})} L({\alpha}^{2} (x)) L(\alpha (y)) L(z)$
\par
$+ (-1)^{\bar{x}\bar{z} + \bar{w}(\bar{x} + \bar{y})} L( \alpha (zx)) L(\alpha (y)) \alpha
- (-1)^{\bar{w}(\bar{x} + \bar{y}) + \bar{y}(\bar{x} + \bar{z})} L(\alpha (y)(xz) {\alpha}^{2}$
\par
$+ (-1)^{\bar{y}(\bar{x} + \bar{w}) + \bar{x}(\bar{z} + \bar{w})} L(\alpha (zy))L(\alpha (x))\alpha$
\par
$- (-1)^{\bar{y}(\bar{x} + \bar{w}) + \bar{x}(\bar{z} + \bar{w})} L({\alpha}^{2} (z)) L(\alpha (y)) L(x) 
= 0$. \\
\\
Switching $x$ and $y$ in (7.5) and next subtracting the result from (7.5), we get 
\par
$(-1)^{\bar{x}(\bar{y} + \bar{w}) + \bar{y}(\bar{w} + \bar{z})} L( {\alpha}^{2} (y)) L(\alpha (x)) L(z) 
+ (-1)^{\bar{y}\bar{z} + \bar{w}(\bar{x} + \bar{y})} L( \alpha (x)(yz)) {\alpha}^{2}$
\par
$+ (-1)^{\bar{w}\bar{y} + \bar{x}(\bar{z} + \bar{w})} L( {\alpha}^{2} (z)) L(\alpha (x)) L(y)
- (-1)^{\bar{w}\bar{x} + \bar{y}(\bar{w} + \bar{z})} L( {\alpha}^{2} (x)) L(\alpha (y)) L(z)$
\par
$-(-1)^{\bar{w}(\bar{x} + \bar{y}) + \bar{y}(\bar{x} + \bar{z})} L(\alpha (y)(xz)) {\alpha}^{2}
- (-1)^{\bar{y}(\bar{x} + \bar{w}) + \bar{x}(\bar{z} + \bar{w})} L( {\alpha}^{2} (z)) L(\alpha (y)) L(x)$ 
\par
$= 0$\\
i.e.
\par
$(-1)^{\bar{w}\bar{x} + \bar{y}(\bar{w} + \bar{z})} [L( {\alpha}^{2} (x)) L(\alpha (y)) L(z)-(-1)^{\bar{x}\bar{y}}
L( {\alpha}^{2} (y)) L(\alpha (x)) L(z)$
\par
$- (-1)^{\bar{z}(\bar{x} + \bar{y})} L( {\alpha}^{2} (z)) L(\alpha (x)) L(y) 
+ (-1)^{\bar{x}\bar{y} + \bar{z}(\bar{x}+\bar{y})} L( {\alpha}^{2} (z)) L(\alpha (y)) L(x)]$
\par
$= (-1)^{\bar{w}\bar{x} + \bar{y}(\bar{w} + \bar{z})} \{ L (\alpha (x)(yz)) {\alpha}^{2}- (-1)^{\bar{x}\bar{y}}
L (\alpha (y)(xz)) {\alpha}^{2} \}$\\
and so, simplifying by $(-1)^{\bar{w}\bar{x} + \bar{y}(\bar{w} + \bar{z})}$ and using (7.4), we obtained
\par
$L(x,y,z) = L(L( \alpha (x))L(y)(z) - (-1)^{\bar{x}\bar{y}} L( \alpha (y))L(x)(z)) {\alpha}^{2}$ \\
which is (7.3). \hfill $\square$ \\
\par
Observe that, in operator language, the Hom-Jordan supertriple identity (see Definition 6.2) looks as\\
\\
(7.6) \; $L_{{\alpha}^{2}(x), {\alpha}^{2}(y)} L_{u,v} - (-1)^{(\bar{u} + \bar{v}) (\bar{x} + \bar{y})} 
L_{{\alpha}^{2}(u), {\alpha}^{2}(v)} L_{x,y} $ 
\par
$- L_{\langle x,y,u \rangle , {\alpha}^{2}(v)} {\alpha}^{2} + (-1)^{\bar{x}\bar{y} + \bar{u} (\bar{x} + \bar{y})}
L_{{\alpha}^{2}(u), \langle y,x,v \rangle} {\alpha}^{2} = 0$. \\
\par
The following lemma computes the first term in (7.6)\\
\par

{\bf Lemma 7.4.} {\it Let $(A, \circ , \alpha)$ be a multiplicative Hom-Jordan superalgebra. Then}\\
\\
(7.7)  $L_{{\alpha}^{2}(x), {\alpha}^{2}(y)} L_{u,v} = L( {\alpha}^{2}(xy)) L ( \alpha (uv)) {\alpha}^{2} + 
L( {\alpha}^{2}(x), {\alpha}^{2}(y)) L(uv) \alpha$
\par
\hspace{3cm}$ + L( {\alpha}^{2}(xy) ) L( \alpha (u), \alpha (v)) \alpha + L( {\alpha}^{2}(x),{\alpha}^{2}(y)) L(u,v)$ \\
\\
{\it for all homogeneous $u,v,x,y \in A$}. \\

\par
{\it Proof.} This lemma is proved in the same way as \cite[Lemma 3.6]{Yau3}. \hfill $\square$ \\
\par
Switching $(x,y)$ with $(u,v)$ in (7.7), we get the expression of the second term in (7.6) as follows. \\
\par

{\bf Lemma 7.5.} {\it Let $(A, \circ , \alpha)$ be a multiplicative Hom-Jordan superalgebra. Then}\\
\\
(7.8) \; $(-1)^{(\bar{u} + \bar{v}) (\bar{x} + \bar{y})} L_{{\alpha}^{2}(u), {\alpha}^{2}(v)} L_{x,y} 
= (-1)^{(\bar{u} + \bar{v}) (\bar{x} + \bar{y})} \{ L( {\alpha}^{2}(uv) ) L( \alpha (xy) ) {\alpha}^{2}$
\par
\hspace{1.0cm}$+ L( {\alpha}^{2}(u), {\alpha}^{2}(v)) L(xy) \alpha + L({\alpha}^{2}(uv))L(\alpha (x),\alpha (y))\alpha$
\par
\hspace{1.0cm}$+ L( {\alpha}^{2}(u), {\alpha}^{2}(v)) L(x,y) \}$ \\
\\
{\it for all homogeneous $u,v,x,y \in A$}. \hfill $\square$ \\

\par
The following lemma computes the sum of the first two terms in (7.6) (i.e. (7.7) and (7.8)). \\
\par

{\bf Lemma 7.6.} {\it Let $(A, \circ , \alpha)$ be a multiplicative Hom-Jordan superalgebra. Then\\
\par
(i) $L( {\alpha}^{2}(xy) )L ( \alpha (uv)) {\alpha}^{2} - (-1)^{(\bar{u} + \bar{v}) (\bar{x} + \bar{y})}
L( {\alpha}^{2}(uv) ) L( \alpha (xy) ) {\alpha}^{2}$
\par
\hspace{0.5cm}$= L( \alpha (xy), \alpha (u)\alpha (v) ) {\alpha}^{2}$;
\par
(ii) $L( {\alpha}^{2}(x), {\alpha}^{2}(y)) L(uv) \alpha - (-1)^{(\bar{u} + \bar{v}) (\bar{x} + \bar{y})}
L({\alpha}^{2}(uv))L(\alpha (x),\alpha (y))\alpha$
\par
\hspace{0.7cm}$= L(L( \alpha (x),\alpha (y) ) (uv)) {\alpha}^{3}$;
\par
(iii) $L( {\alpha}^{2}(xy) ) L( \alpha (u), \alpha (v)) \alpha - (-1)^{(\bar{u} + \bar{v}) (\bar{x} + \bar{y})}
L( {\alpha}^{2}(u), {\alpha}^{2}(v)) L(xy) \alpha$
\par
\hspace{0.7cm}$= (-1)^{(\bar{u} + \bar{v}) (\bar{x} + \bar{y})} \{ (-1)^{\bar{u} \bar{v}} 
L( L({\alpha}^{2}(v))L(\alpha (u))(xy)) {\alpha}^{3}$ 
\par
\hspace{4.0cm}$- L( L({\alpha}^{2}(u))L(\alpha (v))(xy)) {\alpha}^{3} \}$;
\par
(iv) $L( {\alpha}^{2}(x),{\alpha}^{2}(y)) L(u,v) - (-1)^{(\bar{u} + \bar{v}) (\bar{x} + \bar{y})}
L( {\alpha}^{2}(u), {\alpha}^{2}(v)) L(x,y)$
\par
\hspace{0.7cm}$= L(x,y,u,v) - (-1)^{\bar{u}\bar{v} + \bar{x} \bar{y}} L(y,x,v,u)$\\
\\
for all homogeneous $u,v,x,y \in A$}. \\
\par
{\it Proof.} The equality (i) is immediate from the definition of $L(x,y)$. 
\par 
For the equality (ii) observe first that
\par
$L( {\alpha}^{2}(x), {\alpha}^{2}(y)) L(uv)  - (-1)^{(\bar{u} + \bar{v}) (\bar{x} + \bar{y})}
L({\alpha}^{2}(uv))L(\alpha (x),\alpha (y))$ 
\par
$= L( \alpha (x),\alpha (y), uv)$ \\ 
(by the definition of $L(x,y,z)$) and, from the other hand, 
\par
$L(L( \alpha (x),\alpha (y) ) (uv)) {\alpha}^{2} = L( \alpha (x),\alpha (y), uv)$ \\
(by Lemma 7.3) and  so both sides of (ii) are equal to $L(L( \alpha (x),\alpha (y) ) (uv)) \alpha$ which proves (ii).
\par
The equality (iii) holds because both of its sides are equal to \\ $-(-1)^{(\bar{u} + \bar{v}) (\bar{x} + \bar{y})}
L( \alpha (u),\alpha (v) , xy) \alpha$ which is \\ $-(-1)^{(\bar{u} + \bar{v}) (\bar{x} + \bar{y})} L(L(\alpha (u),\alpha (v))(xy))
{\alpha}^{3}$  by Lemma 7.3.
\par
For the equality (iv) observe that, from the definition, 
\par
$L(x,y,u,v) = L( {\alpha}^{2}(x),{\alpha}^{2}(y)) L(\alpha (u))L(v)$ 
\par
\hspace{2.5cm}$- (-1)^{\bar{u} (\bar{x} + \bar{y})}
L({\alpha}^{3}(u)) L( \alpha (x),\alpha (y) )L(v)$
\par
\hspace{2.5cm}$- (-1)^{\bar{v} (\bar{u} + \bar{x} + \bar{y})} L({\alpha}^{3}(v)) L( \alpha (x),\alpha (y) )L(u)$
\par
\hspace{2.5cm}$+ (-1)^{\bar{u}\bar{v} + (\bar{u} + \bar{v}) (\bar{x} + \bar{y})} L({\alpha}^{3}(v))L({\alpha}^{2}(u))L(x,y)$.\\
From the equality above we have
\par
$(-1)^{\bar{u}\bar{v}+\bar{x}\bar{y}} L(y,x,v,u) = (-1)^{\bar{u}\bar{v}+\bar{x}\bar{y}} L({\alpha}^{2}(y),{\alpha}^{2}(x))L(\alpha (v))L(u)$
\par
\hspace{4.2cm}$- (-1)^{\bar{u}\bar{v}+\bar{x}\bar{y} + \bar{v} (\bar{x} + \bar{y})} 
L({\alpha}^{3}(v)) L( \alpha (y),\alpha (x) )L(u)$
\par
\hspace{4.2cm}$- (-1)^{\bar{x}\bar{y}+\bar{u} (\bar{x} + \bar{y})} L({\alpha}^{3}(u)) L(\alpha (y),\alpha (x))L(v)$
\par
\hspace{4.2cm}$+ (-1)^{\bar{x}\bar{y}+ (\bar{u}+\bar{v})(\bar{x}+\bar{y})}L({\alpha}^{3}(u))L({\alpha}^{2}(v))L(y,x)$.\\
Therefore, using Lemma 7.1, we obtain
\par
$L(x,y,u,v) + (-1)^{\bar{u}\bar{v}+\bar{x}\bar{y}} L(y,x,v,u)$
\par
$= L({\alpha}^{2}(x),{\alpha}^{2}(y))
\{ L(\alpha (u))L(v) - (-1)^{\bar{u}\bar{v}} L(\alpha (v))L(u) \}$
\par
\hspace{0.4cm}$- (-1)^{(\bar{u} + \bar{v}) (\bar{x} + \bar{y})} 
\{ L({\alpha}^{3}(u))L({\alpha}^{2}(v)) - (-1)^{\bar{u}\bar{v}} L({\alpha}^{3}(v))L({\alpha}^{2}(u)) \} L(x,y)$\\
which proves the equality (iv) by Definition 7.1. \hfill $\square$ \\
\par
The following lemma computes the third term in (7.6). \\
\par

{\bf Lemma 7.7.} {\it Let $(A, \circ , \alpha)$ be a multiplicative Hom-Jordan superalgebra. Then}\\
\par
$L_{\langle x,y,u \rangle , {\alpha}^{2}(v)} {\alpha}^{2} = (-1)^{\bar{u}\bar{v}+(\bar{u}+\bar{v}) (\bar{x}+\bar{y})}
L( L({\alpha}^{2}(v)) L(\alpha (u))(xy)) {\alpha}^{3}$
\par
\hspace{3cm}$+ L(L(\alpha (x), \alpha (y)) (uv)) {\alpha}^{3}$
\par
\hspace{3cm}$- (-1)^{\bar{u} (\bar{x} + \bar{y})} L( L({\alpha}^{2}(u)) L(x,y)(v)) {\alpha}^{3}$
\par
\hspace{3cm}$- (-1)^{\bar{v} (\bar{u} + \bar{x} + \bar{y})} L({\alpha}^{2}(v), (xy) \alpha (u)) {\alpha}^{2} + L(x,y,u,v)$\\
\\
{\it for all homogeneous $u,v,x,y \in A$}. \\

\par
{\it Proof.} We have
\par
$L_{\langle x,y,u \rangle , {\alpha}^{2}(v)} {\alpha}^{2} = L( \langle x,y,u \rangle {\alpha}^{2}(v)) {\alpha}^{3}
+L( \langle x,y,u \rangle , {\alpha}^{2}(v)) {\alpha}^{2} $ (by Definition 7.1)
\par
\hspace{2.5cm}$=L(L(L(xy)(\alpha (u)))( {\alpha}^{2}(v) )) {\alpha}^{3} + L(L(L(x,y)(u)) ({\alpha}^{2}(v))) {\alpha}^{3}$
\par
\hspace{2.5cm}$+ L(\alpha (\langle x,y,u \rangle)) L({\alpha}^{2}(v)) {\alpha}^{2}$ 
\par
\hspace{2.5cm}$- (-1)^{\bar{v} (\bar{u} + \bar{x} + \bar{y})}
L({\alpha}^{3}(v))  L( \langle x,y,u \rangle ) {\alpha}^{2}$
\par
\hspace{2.5cm}$=(-1)^{\bar{v} (\bar{u} + \bar{x} + \bar{y}) + \bar{u}(\bar{x} + \bar{y})} 
L(L({\alpha}^{2}(v)) L(\alpha (u))(xy)) {\alpha}^{3}$ 
\par
\hspace{2.5cm}$+ L(L(x,y,u)(v)) {\alpha}^{3}$
\par
\hspace{2.5cm}$+L(\alpha (L(xy)(\alpha (u)) + L(x,y)(u)))L({\alpha}^{2}(v)) {\alpha}^{2}$
\par
\hspace{2.5cm}$- (-1)^{\bar{v} (\bar{u} + \bar{x} + \bar{y})} L({\alpha}^{3}(v)) L(L(xy)(\alpha (u)) + L(x,y)(u)){\alpha}^{2}$
(by supercommutativity, (7.2) and Lemma 7.3)
\par
\hspace{2.5cm}$= (-1)^{\bar{u}\bar{v} + (\bar{u} + \bar{v})(\bar{x} + \bar{y})} 
L(L({\alpha}^{2}(v)) L(\alpha (u))(xy)) {\alpha}^{3}$ 
\par
\hspace{2.5cm}$+ L(L(x,y,u)(v)) {\alpha}^{3}$
\par
\hspace{2.5cm}$+ L(\alpha (xy) {\alpha}^{2}(u)) L({\alpha}^{2}(v)) {\alpha}^{2} + L( \alpha (x), \alpha (y), \alpha (u))L(v)$
\par
\hspace{2.5cm}$- (-1)^{\bar{v} (\bar{u} + \bar{x} + \bar{y})} L({\alpha}^{3}(v)) L(L(xy)(\alpha (u))) {\alpha}^{2}$
\par
\hspace{2.5cm}$- (-1)^{\bar{v} (\bar{u} + \bar{x} + \bar{y})} L({\alpha}^{3}(v))L(x,y,u)$ (using repeatedly Lemmas 7.1 and 7.3)
\par
\hspace{2.5cm}$= (-1)^{\bar{u}\bar{v} + (\bar{u} + \bar{v})(\bar{x} + \bar{y})}
L(L({\alpha}^{2}(v)) L(\alpha (u))(xy)) {\alpha}^{3}$ 
\par
\hspace{2.5cm}$+ L((L( \alpha (x), \alpha (y))L(u)$
\par
\hspace{2.5cm}$- (-1)^{\bar{u} (\bar{x} + \bar{y})} L({\alpha}^{2}(u)) L(x,y)(v)) {\alpha}^{3}$
\par
\hspace{2.5cm}$- (-1)^{\bar{v} (\bar{u} + \bar{x} + \bar{y})} \{ L({\alpha}^{3}(v)) L(L(xy)(\alpha (u)))$ 
\par
\hspace{5cm}$- (-1)^{\bar{v} (\bar{u} + \bar{x} + \bar{y})} L( \alpha (xy) {\alpha}^{2}(u)) L({\alpha}^{2}(v))\} {\alpha}^{2}$
\par
\hspace{2.5cm}$+L( \alpha (x), \alpha (y), \alpha (u))L(v) - (-1)^{\bar{v} (\bar{u} + \bar{x} + \bar{y})}
L({\alpha}^{3}(v))L(x,y,u)$ (rearranging terms and using the definiton of $L(x,y,z)$)
\par
\hspace{2.5cm}$=(-1)^{\bar{u}\bar{v} + (\bar{u} + \bar{v})(\bar{x} + \bar{y})}
L(L({\alpha}^{2}(v)) L(\alpha (u))(xy)) {\alpha}^{3}$
\par
\hspace{2.5cm}$+ L(L( \alpha (x), \alpha (y))(uv)) {\alpha}^{3}$
\par
\hspace{2.5cm}$-(-1)^{\bar{u} (\bar{x} + \bar{y})} L( L({\alpha}^{2}(u)) L(x,y)(v)) {\alpha}^{3}$
\par
\hspace{2.5cm}$- (-1)^{\bar{v} (\bar{u} + \bar{x} + \bar{y})} L({\alpha}^{2}(v), (xy) \alpha (u)) {\alpha}^{2} + L(x,y,u,v)$ (by
Definition 7.1) \\
and this proves the lemma. \hfill $\square$ \\
\par
As in the proof of Lemma 7.7, we compute the fourth term in (7.6) as follows.\\
\par

{\bf Lemma 7.8.} {\it Let $(A, \circ , \alpha)$ be a multiplicative Hom-Jordan superalgebra. Then}\\
\par
$(-1)^{\bar{x}\bar{y} + \bar{u} (\bar{x} + \bar{y})}
L_{{\alpha}^{2}(u), \langle y,x,v \rangle} {\alpha}^{2}= (-1)^{(\bar{u}+\bar{v}) (\bar{x}+\bar{y})}
L(L({\alpha}^{2}(u)) L(\alpha (v))(xy)) {\alpha}^{3}$
\par
\hspace{5.3cm}$- (-1)^{\bar{u} (\bar{x} + \bar{y})} L(L({\alpha}^{2}(u))L(x,y)(v)) {\alpha}^{3}$
\par
\hspace{5.3cm}$+ (-1)^{(\bar{u}+\bar{v}) (\bar{x}+\bar{y})} L({\alpha}^{2}(u), \alpha (v) (xy)) {\alpha}^{2}$
\par
\hspace{5.3cm}$- (-1)^{\bar{u}\bar{v} + \bar{x}\bar{y}} L(y,x,v,u)$\\
\\
{\it for all homogeneous $u,v,x,y \in A$}. \\

\par
{\it Proof.} We have
\par
$L_{{\alpha}^{2}(u), \langle y,x,v \rangle} {\alpha}^{2}= L(L({\alpha}^{2}(u))L(yx)(\alpha (v))) {\alpha}^{3}
+ L(L({\alpha}^{2}(u))L(y,x)(v)) {\alpha}^{3}$
\par
\hspace{2.5cm}$+L({\alpha}^{3}(u)) L(\langle y,x,v \rangle ) {\alpha}^{2}$
\par
\hspace{2.5cm}$- (-1)^{\bar{u} (\bar{v} + \bar{x} + \bar{y})} L(\alpha (\langle y,x,v \rangle )) L({\alpha}^{2}(u)) {\alpha}^{2}$
\par
\hspace{2.5cm}$=L(L({\alpha}^{2}(u))L(yx)(\alpha (v))) {\alpha}^{3}
+ L(L({\alpha}^{2}(u))L(y,x)(v)) {\alpha}^{3} $
\par
\hspace{2.5cm}$+ L({\alpha}^{3}(u)) L(L(yx)(\alpha (v))) {\alpha}^{2} + L(L({\alpha}^{3}(u))L(y,x)(v)) {\alpha}^{2}$
\par
\hspace{2.5cm}$- (-1)^{\bar{u} (\bar{v} + \bar{x} + \bar{y})} L(\alpha (L(yx)(\alpha (v))) L({\alpha}^{2}(u)) {\alpha}^{2}$
\par
\hspace{2.5cm}$- (-1)^{\bar{u} (\bar{v} + \bar{x} + \bar{y})} L(\alpha (L(y,x)(v)))  L({\alpha}^{2}(u)) {\alpha}^{2}$
\par
\hspace{2.5cm}$=L(L({\alpha}^{2}(u))L(yx)(\alpha (v))) {\alpha}^{3}
+ L(L({\alpha}^{2}(u))L(y,x)(v)) {\alpha}^{3}$
\par
\hspace{2.5cm}$+ L({\alpha}^{3}(u)) L((yx)(\alpha (v))) {\alpha}^{2}$ 
\par
\hspace{2.5cm}$- (-1)^{\bar{x}\bar{y}+\bar{u}(\bar{v}+\bar{x}+\bar{y})} 
L( \alpha (xy) {\alpha}^{2}(v)) L({\alpha}^{2}(u)) {\alpha}^{2}$ 
\par
\hspace{2.5cm}$+ L({\alpha}^{3}(u)) L(y,x,v)$
\par
\hspace{2.5cm}$- (-1)^{\bar{u} (\bar{v} + \bar{x} + \bar{y})} 
L(\alpha (y),\alpha (x),\alpha (v))L(u)$ (by Lemma 7.3)
\par
\hspace{2.5cm}$= (-1)^{\bar{x}\bar{y}+\bar{v}(\bar{x}+\bar{y})} L(L({\alpha}^{2}(u))L(\alpha (v))(xy)){\alpha}^{3}$
\par
\hspace{2.5cm}$- (-1)^{\bar{x}\bar{y}} L(L({\alpha}^{2}(u))L(x,y)(v)) {\alpha}^{3}$
\par
\hspace{2.5cm}$+ (-1)^{\bar{x}\bar{y}+\bar{v}(\bar{x}+\bar{y})} \{ L({\alpha}^{3}(u))L(\alpha (v))(xy))$
\par
\hspace{5.5cm}$- (-1)^{\bar{u} (\bar{v} + \bar{x} + \bar{y})}L({\alpha}^{2}(v)\alpha (xy))L({\alpha}^{2}(u)) \} {\alpha}^{2}$
\par
\hspace{2.5cm}$- (-1)^{\bar{u} (\bar{v} + \bar{x} + \bar{y})} \{ L(\alpha (y), \alpha (x), \alpha (v))L(u)$
\par
\hspace{5.5cm}$- (-1)^{\bar{u} (\bar{v} + \bar{x} + \bar{y})} L({\alpha}^{3}(u)) L(y,x,v) \}$\\
and the desired equality follows from Definition 7.1. \hfill $\square$ \\
\par
We are now in position to prove Theorem 7.1.\\
\par
{\it Proof of Theorem 7.1.} According to Definition 6.2 we need to verify the outer symmetry and the Hom-Jordan 
supertriple identity (7.6) for \\ $(A, \langle \cdot , \cdot , \cdot \rangle, {\alpha}^{2})$.
\par
The outer symmetry clearly follows from (7.1) and the supercommutativity of $(A, \circ , \alpha)$. It remains to prove
(7.6) for $(A, \langle \cdot , \cdot , \cdot \rangle, {\alpha}^{2})$.
\par
Using Lemma 7.4 - Lemma 7.8, we find that the left hand side of (7.6) reduces to
\par
$\{ L(\alpha (xy), \alpha (u) \alpha (v)) + (-1)^{(\bar{u} + \bar{v})(\bar{x} + \bar{y})} L({\alpha}^{2}(u),\alpha (v)(xy))$
\par
\hspace{5.5cm}$+ (-1)^{\bar{v} (\bar{u} + \bar{x} + \bar{y})} L({\alpha}^{2}(v), (xy) \alpha (u)) \} {\alpha}^{2}$ \\
i.e.
\par
${\circlearrowleft}_{\alpha (u),\alpha (v),xy} (-1)^{\bar{u} (\bar{v} + \bar{x} + \bar{y})} L( \alpha (xy), \alpha (u) \alpha (v)) {\alpha}^{2}$ \\
which is equal to $0$ by Lemma 7.2. It then follows that $(A, \langle \cdot , \cdot , \cdot \rangle, {\alpha}^{2})$
is a multiplicative Hom-Jordan supertriple system. \hfill $\square$

\subsection{From right Hom-alternative superalgebras to Hom-Bol superalgebras}

\par
In this subsection we prove a $\mathbb{Z}_2$-graded generalization of results connecting right Hom-alternative algebras and 
Hom-Bol algebras \cite{AI2}. Specifically we aim to prove that every multiplicative right Hom-alternative superalgebra
has a natural Hom-Bol structure, extending thusly Theorem 3.1. First we need the following result which is a 
$\mathbb{Z}_2$-graded generalization of \cite[Theorem 4.2]{Yau5}.\\
\par

{\bf Theorem 7.2.} {\it Let $(V, \langle \cdot , \cdot , \cdot \rangle, \theta)$ be any multiplicative Hom-Jordan
supertriple system. If define on $V$ a ternary operation $[x,y,z] := \langle x,y,z \rangle - (-1)^{\bar{x}\bar{y}}
\langle y,x,z \rangle$ for all homogeneous $x,y,z \in V$, then $(V, [ \cdot , \cdot , \cdot ], \theta)$ is a 
multiplicative Hom-Lie supertriple system}.\\

\par
{\it Proof.} First observe that the multiplicativity of $(V, [ \cdot , \cdot , \cdot ], \theta)$ follows from the one
of $(V, \langle \cdot , \cdot , \cdot \rangle, \theta)$  and the checking of (SHB2) and (SHB3) is straightforward.
For the checking of the identity (SHB5), note that its left-hand side is
\par
$[\theta (x), \theta (y), [u,v,w]] = \langle \theta (x), \theta (y), [u,v,w] \rangle - (-1)^{\bar{x}\bar{y}}
\langle \theta (y), \theta (x), [u,v,w] \rangle$
\par
\hspace{0.5cm}$= \langle \theta (x), \theta (y), \langle u,v,w \rangle \rangle - (-1)^{\bar{u}\bar{v}}
\langle \theta (x), \theta (y), \langle v,u,w \rangle \rangle$
\par
\hspace{0.5cm}$- (-1)^{\bar{x}\bar{y}} \langle \theta (y), \theta (x), \langle u,v,w \rangle \rangle
+ (-1)^{\bar{u}\bar{v}+\bar{x}\bar{y}} \langle \theta (y), \theta (x), \langle v,u,w \rangle \rangle$
\par
\hspace{0.5cm}$=\langle \langle x,y,u \rangle, \theta (v), \theta (w) \rangle + (-1)^{(\bar{u} + \bar{v})(\bar{x} + \bar{y})}
(\langle \theta (u), \theta (v), \langle x,y,w \rangle \rangle$
\par
\hspace{0.5cm}$- \langle \theta (u), \langle v,x,y \rangle, \theta (w) \rangle) - (-1)^{\bar{u}\bar{v}}
\langle \langle x,y,v \rangle, \theta (u), \theta (w) \rangle$
\par
\hspace{0.5cm}$- (-1)^{\bar{u}\bar{v}+ (\bar{u} + \bar{v})(\bar{x} + \bar{y})} (\langle \theta (v), \theta (u), \langle x,y,w \rangle \rangle
- \langle \theta (v), \langle u,x,y \rangle, \theta (w) \rangle)$
\par
\hspace{0.5cm}$- (-1)^{\bar{x}\bar{y}} \langle \langle y,x,u \rangle, \theta (v), \theta (w) \rangle$
\par
\hspace{0.5cm}$- (-1)^{\bar{x}\bar{y}+ (\bar{u} + \bar{v})(\bar{x} + \bar{y})}
(\langle \theta (u), \theta (v), \langle y,x,w \rangle \rangle 
- \langle \theta (u),\langle v,y,x \rangle,\theta (w) \rangle)$
\par
\hspace{0.5cm}$+ (-1)^{\bar{u}\bar{v}+\bar{x}\bar{y}} 
\langle \langle y,x,v \rangle, \theta (u), \theta (w) \rangle$
\par
\hspace{0.5cm}$+ (-1)^{\bar{u}\bar{v}+\bar{x}\bar{y} + (\bar{u} + \bar{v})(\bar{x} + \bar{y})}
(\langle \theta (v), \theta (u), \langle y,x,w \rangle \rangle - \langle \theta (v),\langle u,y,x \rangle,\theta (w) \rangle)$\\
and so \\
\\
(7.9) \; $[\theta (x), \theta (y), [u,v,w]] =$ 
\par
\hspace{0.5cm}$\langle \langle x,y,u \rangle, \theta (v), \theta (w) \rangle
+ (-1)^{(\bar{u} + \bar{v})(\bar{x} + \bar{y})} \langle \theta (u), \theta (v), \langle x,y,w \rangle \rangle$
\par
\hspace{0.5cm}$- (-1)^{(\bar{u} + \bar{v})(\bar{x} + \bar{y})} \langle \theta (u), \langle v,x,y \rangle, \theta (w) \rangle
- (-1)^{\bar{u}\bar{v}} \langle \langle x,y,v \rangle, \theta (u), \theta (w) \rangle$
\par
\hspace{0.5cm}$- (-1)^{\bar{u}\bar{v}+ (\bar{u} + \bar{v})(\bar{x} + \bar{y})} 
\langle \theta (v), \theta (u), \langle x,y,w \rangle \rangle$ 
\par
\hspace{0.5cm}$+ (-1)^{\bar{u}\bar{v}+ (\bar{u} + \bar{v})(\bar{x} + \bar{y})}
\langle \theta (v), \langle u,x,y \rangle, \theta (w) \rangle$
\par
\hspace{0.5cm}$-(-1)^{\bar{x}\bar{y}} \langle \langle y,x,u \rangle, \theta (v), \theta (w) \rangle
- (-1)^{\bar{x}\bar{y}+ (\bar{u} + \bar{v})(\bar{x} + \bar{y})} \langle \theta (u), \theta (v), \langle y,x,w \rangle \rangle$
\par
\hspace{0.5cm}$+(-1)^{\bar{x}\bar{y}+ (\bar{u} + \bar{v})(\bar{x} + \bar{y})}
\langle \theta (u),\langle v,y,x \rangle,\theta (w) \rangle + (-1)^{\bar{u}\bar{v}+\bar{x}\bar{y}} 
\langle \langle y,x,v \rangle, \theta (u), \theta (w) \rangle$
\par
\hspace{0.5cm}$+(-1)^{\bar{u}\bar{v}+\bar{x}\bar{y} + (\bar{u} + \bar{v})(\bar{x} + \bar{y})}
\langle \theta (v), \theta (u), \langle y,x,w \rangle \rangle$
\par
\hspace{0.5cm}$- (-1)^{\bar{u}\bar{v}+\bar{x}\bar{y} + (\bar{u} + \bar{v})(\bar{x} + \bar{y})}
\langle \theta (v),\langle u,y,x \rangle,\theta (w) \rangle$. \\
\\
Likewise we compute
\par
$[[x,y,u], \theta (v), \theta (w)] = \langle \langle x,y,u \rangle, \theta (v), \theta (w) \rangle
- (-1)^{\bar{x}\bar{y}} \langle \langle y,x,u \rangle, \theta (v), \theta (w) \rangle$
\par
\hspace{0.5cm}$-(-1)^{\bar{v} (\bar{u}+\bar{x}+\bar{y})} \langle \theta (v),\langle x,y,u \rangle,\theta (w) \rangle
+(-1)^{\bar{x}\bar{y}+\bar{v}(\bar{u}+\bar{x}+\bar{y})} \langle \theta (v),\langle y,x,u \rangle,\theta (w) \rangle$;
\par
$(-1)^{\bar{u} (\bar{x} + \bar{y})} [ \theta (u), [x,y,v], \theta (w)] = (-1)^{\bar{u} (\bar{x} + \bar{y})}
\langle \theta (u),\langle x,y,v \rangle,\theta (w) \rangle$
\par
\hspace{5.5cm}$- (-1)^{\bar{x}\bar{y}+ \bar{u} (\bar{x} + \bar{y})}
\langle \theta (u),\langle y,x,v \rangle,\theta (w) \rangle$
\par
\hspace{5.5cm}$-(-1)^{\bar{u}\bar{v}} \langle \langle x,y,v \rangle, \theta (u), \theta (w) \rangle$
\par
\hspace{5.5cm}$+ (-1)^{\bar{u}\bar{v}+ \bar{x}\bar{y}} \langle \langle y,x,v \rangle,\theta (u),\theta (w) \rangle$\\
and \\
$(-1)^{(\bar{u} + \bar{v})(\bar{x} + \bar{y})} [\theta (u), \theta (v), [x,y,w]]
= (-1)^{(\bar{u} + \bar{v})(\bar{x} + \bar{y})} \langle \theta (u), \theta (v), \langle x,y,w \rangle \rangle$
\par
\hspace{5.5cm}$- (-1)^{\bar{x}\bar{y}+ (\bar{u}+ \bar{v})(\bar{x} + \bar{y})} 
\langle \theta (u), \theta (v), \langle y,x,w \rangle \rangle$
\par
\hspace{5.5cm}$- (-1)^{\bar{u}\bar{v}+ (\bar{u}+ \bar{v})(\bar{x} + \bar{y})} 
\langle \theta (v), \theta (u), \langle x,y,w \rangle \rangle$
\par
\hspace{5.5cm}$+ (-1)^{\bar{u}\bar{v} + \bar{x}\bar{y} + (\bar{u}+ \bar{v})(\bar{x} + \bar{y})} 
\langle \theta (v), \theta (u), \langle y,x,w \rangle \rangle$. \\
Therefore, summing up the expressions of $[[x,y,u], \theta (v), \theta (w)]$, \\ 
$(-1)^{\bar{u} (\bar{x} + \bar{y})} [ \theta (u), [x,y,v], \theta (w)]]$,
and $(-1)^{(\bar{u} + \bar{v})(\bar{x} + \bar{y})} [\theta (u), \theta (v), [x,y,w]]$ as above and using the outer
supersymmetry whenever applicable, we obtain the twelve terms of the right-hand side of (7.9). So (SHB5) is verified
and thus $(V, \langle \cdot , \cdot , \cdot \rangle, \theta)$ is a Hom-Lie supertriple system. \hfill $\square$ \\
\par
We can now prove the main result of this section. \\
\par
{\bf Theorem 7.3.} {\it The supercommutator Hom-superalgebra of any multiplicative right Hom-alternative superalgebra is a
Hom-Bol superalgebra}.\\
\par
{\it Proof.} Let ${\mathcal A} := (A,*, \alpha)$ be a multiplicative right Hom-alternative superalgebra. Then, by
Theorem 4.1, ${\mathcal A}^{+}$ is a multiplicative Hom-Jordan superalgebra and Theorem 7.1 implies that
$(A, \langle \cdot, \cdot , \cdot \rangle, {\alpha}^{2})$ is a Hom-Jordan supertriple system with $\langle \cdot, \cdot , \cdot \rangle$
defined by (7.1). Next Theorem 7.2 says that $(A, [ \cdot, \cdot , \cdot ], {\alpha}^{2})$ is a Hom-Lie supertriple
system. Now define on $A$ the ternary product
\par
$\{ x,y,z \} := {\frac{1}{2}} [x,y,z]$. \\
Then $(A, \{ \cdot, \cdot , \cdot \}, {\alpha}^{2})$ is also a Hom-Lie supertriple system and so (SHB2), (SHB3), and
(SHB5) hold. To conclude, we only need to prove that (SHB4) holds since (SHB01), (SHB02 and (SHB1) are quite obvious.
\par
First observe that $\{ x,y,z \} = (-1)^{\bar{x}(\bar{y} + \bar{z})} as_{{\mathcal A}^{+}} (y,z,x)$ and it is easily
seen that the Hom-version of (3.2) is \\
\\
(7.10) \; $[x,y,z] = 2[[x,y], \alpha (z)] - (-1)^{\bar{z}(\bar{x} + \bar{y})} as_{\mathcal A} (z,x,y)$, \\
\\
where $as_{\mathcal A} (z,x,y) = (z*x)*\alpha (y) - \alpha (z) * (x*y)$ (a proof of (7.10) is similar to that of (3.2)
in Lemma 3.1). From (7.10) we get \\
\\
(7.11) \; ${\frac{1}{2}}as_{\mathcal A} (z,x,y) = (-1)^{\bar{z}(\bar{x} + \bar{y})} ( [[x,y], \alpha (z)] - \{ x,y,z \} )$\\
\\
(see (3.4) for the untwisted version). Now, using (7.11), (2.6) and proceeding as in the proof of Theorem 3.1 we
obtain (SHB4). Thus $(A, [ \cdot, \cdot ], \{ \cdot, \cdot , \cdot \},{\alpha}^{2})$ is a Hom-Bol superalgebra. This
completes the proof. \hfill $\square$
\\

\end{document}